\newif\ifJOC

\JOCfalse 

\documentclass[11pt, 3p]{elsarticle}

\makeatletter
\def\ps@pprintTitle{%
 \let\@oddhead\@empty
 \let\@evenhead\@empty
 \def\@oddfoot{}%
 \let\@evenfoot\@oddfoot}
\makeatother

\usepackage{lineno}
\usepackage{amstext}
\usepackage{amsmath,amssymb,amsthm,graphicx}
\modulolinenumbers[5]

\DeclareMathOperator*{\argmax}{arg\,max}

\newtheorem{theorem}{Theorem}[section]
\newtheorem{definition}[theorem]{Definition}
\newtheorem{lemma}[theorem]{Lemma}

\newtheorem{proposition}[theorem]{Proposition}

\theoremstyle{remark}

\usepackage{bm}
\usepackage{amsfonts}
\usepackage[short]{optidef}
\usepackage{changepage}
\usepackage[utf8x]{inputenc}
\usepackage[T1]{fontenc}
\usepackage{lmodern}
\usepackage[english]{babel}
\usepackage{hyperref}
\usepackage{comment}
\usepackage{longtable}
\usepackage{enumitem}
\usepackage{algpseudocode}
\usepackage[ruled,vlined]{algorithm2e}
\usepackage{subcaption}
\usepackage{pgfplots}
\usepackage{booktabs}
\usepackage{tabularx}
\usepackage{ltablex}
\usepackage{longtable}
\usepackage{comment}
\usepackage{url}
\pgfplotsset{compat=1.9}
\usepackage{nicematrix}
\usepackage{bbm}

\newcommand{\inprod}[2]{{\left \langle #1,#2 \right \rangle}} 

\newcommand{\Zbar}{\bar{Z}}

\newcommand{\st}{\textrm{s.t.}}
\newcommand{\tr}{\textrm{tr}}
\newcommand{\Diag}{\textrm{Diag}}
\newcommand{\diag}{\textrm{diag}}

\newcommand{\rank}{\textrm{rank}}

\biboptions{authoryear}

\usepackage{setspace}
\usepackage{background}

\begin{document}


\SetBgContents{Published at \url{https://doi.org/10.1287/ijoc.2024.0683}}      
\SetBgPosition{current page.center}
\SetBgAngle{0}                                    
\SetBgColor{gray}                                 
\SetBgScale{1.5}                                  
\SetBgHshift{0}                                   
\SetBgVshift{9cm} 

\begin{frontmatter}
\title{An SDP-based Branch-and-Cut Algorithm for Biclustering}


\author{Antonio M. Sudoso}
\ead{antoniomaria.sudoso@uniroma1.it}
\address{Department of Computer, Control and Management Engineering ``Antonio Ruberti'', \\ Sapienza University of Rome, Via Ariosto 25, 00185, Italy}
\begin{abstract}
Biclustering, also called co-clustering, block clustering, or two-way clustering, involves the simultaneous clustering of both the rows and columns of a data matrix into distinct groups, such that the rows and columns within a group display similar patterns.
As a model problem for biclustering, we consider the $k$-densest-disjoint biclique problem, whose goal is to identify $k$ disjoint complete bipartite subgraphs (called bicliques) of a given weighted complete bipartite graph such that the sum of their densities is maximized.
To address this problem, we present a tailored branch-and-cut algorithm. For the upper bound routine, we consider a semidefinite programming relaxation and propose valid inequalities to strengthen the bound. We solve this relaxation in a cutting-plane fashion using a first-order method. For the lower bound, we design a maximum weight matching rounding procedure that exploits the solution of the relaxation solved at each node.
Computational results on both synthetic and real-world instances show that the proposed algorithm can solve instances approximately 20 times larger than those handled by general-purpose solvers.
\end{abstract}


\end{frontmatter}

\section{Introduction}
Cluster analysis is a fundamental task in data science, optimization, and computer science. Its goal is to partition a set of entities 
into so-called clusters such that entities within each cluster are more similar to each other than to those in other clusters.
In many applications, there is often a need to cluster the entries of the data matrix into different groups such that the rows (samples) and columns (features) within a group exhibit similar patterns. This task is recognized in the literature by various names, including biclustering, block clustering, co-clustering, two-way clustering, or two-mode clustering \citep{busygin2008biclustering}.
The term \textit{biclustering} was initially introduced by \cite{mirkin1996mathematical} and notably employed by \cite{cheng2000biclustering} in the analysis of gene expression data, playing a significant role in popularizing biclustering techniques. However, the roots of biclustering can be traced back to the work of \cite{hartigan1972direct}.


In {contrast} to one-way clustering methods (e.g., hierarchical or $k$-means clustering), which can only locate global patterns, biclustering can also discover local patterns. Traditional clustering applied to the rows (samples) of the data matrix considers all available columns (features). On the other hand, biclustering involves clustering rows into different groups based on various subsets of columns. This approach allows the identification of subsets of samples exhibiting similar behavior across a specific set of features, and vice versa \citep{madeira2004biclustering}.
This versatility makes biclustering particularly well-suited for applications such as market segmentation, gene expression analysis, and document clustering. For instance, biclustering can be employed to identify specific subregions within 2D gene expression data, where rows represent genes and columns represent either tissue samples or experimental conditions \citep{cheng2000biclustering}.
In document clustering, the objective is to identify patterns within the content of documents that may not be apparent through traditional clustering methods. This involves clustering both documents and terms simultaneously, revealing relationships between certain topics and terms across a subset of documents \citep{dhillon2001co}.
In market segmentation, the focus is on detecting subgroups of customers with similar preferences toward corresponding subgroups of products \citep{hofmann1999collaborative}. For an extensive overview of biclustering applications, one can refer to the survey conducted by \cite{busygin2008biclustering} and \cite{fan2010recent}.

\begin{figure}[!ht]
\centering
\begin{subfigure}{.3\textwidth}
\centering
\includegraphics[width=0.75\linewidth]{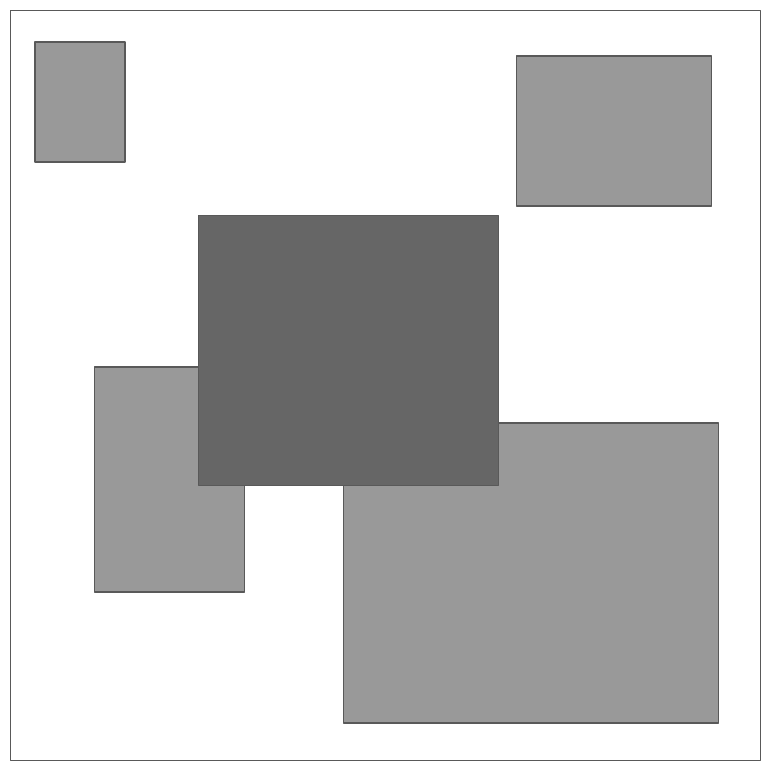}
\caption{}
\label{overlapping}
\end{subfigure}
\begin{subfigure}{.3\textwidth}
\centering
\includegraphics[width=0.75\linewidth]{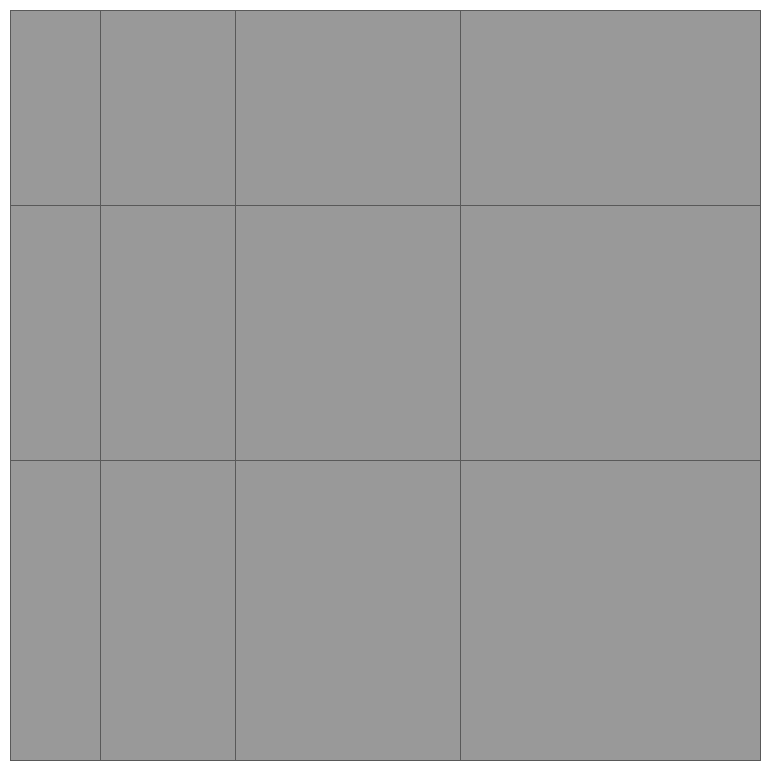}
\caption{}
\label{checkerboard}
\end{subfigure}
\begin{subfigure}{.3\textwidth}
\centering
\includegraphics[width=0.75\linewidth]{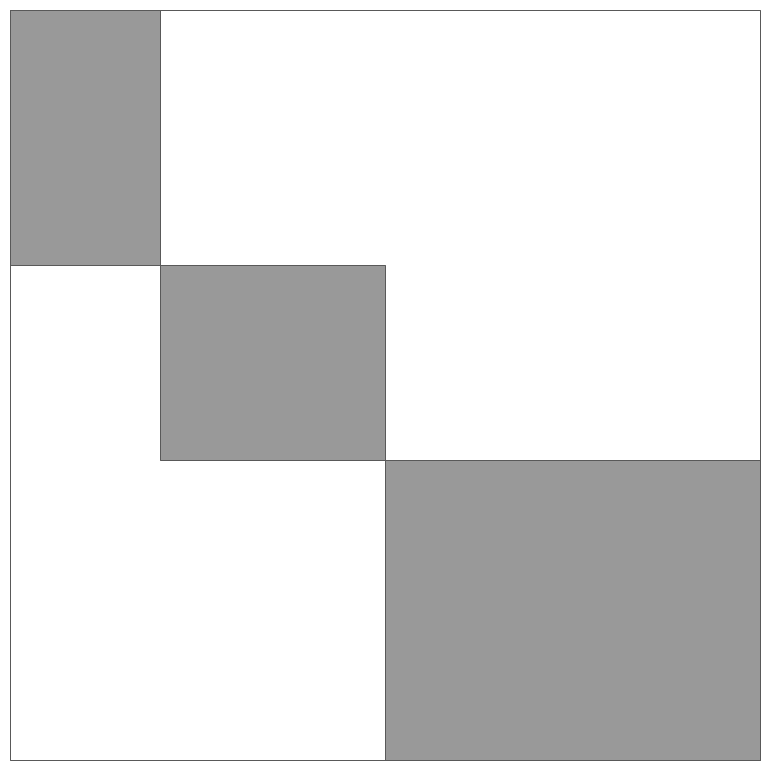}
\caption{}
\label{exclusive}
\end{subfigure}
\caption{Different types of bicluster structures (after row and column reordering): (a) Overlapping biclusters, (b) non-overlapping biclusters with checkerboard structure, and (c) exclusive rows and columns biclusters with block-diagonal structure. The algorithm proposed in this paper identifies bicluster structures of type (c).}
\label{fig:bcstructure}
\end{figure}

Over the years, a large number of biclustering methods have been proposed, each one characterized by different attributes, such as computational complexity, interpretability, and optimization
criterion. Comprehensive reviews of these methods can be found in \citep{madeira2004biclustering} and \citep{padilha2017systematic}. Although biclustering methods differ considerably, they can be classified into different groups based on the clustering structure they are able to reveal. Figure \ref{fig:bcstructure} shows three types of bicluster structures. Specifically, when biclustering algorithms assume the existence of several blocks in the data matrix, the following structures can be obtained: (i) overlapping biclusters, (ii) non-overlapping biclusters with checkerboard structure, and (iii) exclusive rows and columns biclusters with block-diagonal structure. 
In this work, we aim to identify biclusters where a row and a column can only belong to a single bicluster, resulting in exclusive row and column biclusters. This approach ensures straightforward interpretation and visualization. Specifically, rearranging the rows and columns of the data matrix reveals the biclusters on the diagonal, a pattern that naturally emerges in various applications \citep{busygin2008biclustering}.
For instance, in gene expression analysis, block-diagonal biclusters simplify the categorization of genes into different groups associated with diverse types of cancer tissues. Similarly, in text mining, block-diagonal biclusters aid in grouping documents into distinct clusters based on unique sets of terms \citep{dhillon2001co}.

\section{Related Work}
Biclustering is closely connected to mathematical optimization and graph theory. 
The complexity of finding an optimal biclustering depends on the measure of fitness of a proposed partition. The most interesting models for biclustering are posed as an intractable combinatorial optimization problem \citep{ding2006biclustering, fan2012multi, ames2014guaranteed}. However, owing to its practical importance, several heuristic methods using local search strategies have been developed. For instance, a popular approach is to search for one bicluster at a time by iteratively assigning rows to a bicluster based on
the columns, and vice versa \citep{bergmann2003iterative, shabalin2009finding}. 

{Many biclustering methods are related to matrix decomposition techniques such as the Singular Value Decomposition (SVD) \citep{busygin2008biclustering}. Although not specifically designed for biclustering, these methods can identify biclusters through sparsity constraints, viewing biclustering as sparse matrix factorization. Sparse truncated SVD (SSVD), introduced by \cite{witten2009penalized}, has been used by \cite{lee2010biclustering} to identify interpretable row–column associations within high-dimensional data matrices. SSVD seeks a low-rank, checkerboard-structured matrix approximation by making both the left- and right-singular vectors orthonormal and sparse. This is achieved by treating the singular vectors of a regular SVD as coefficients of a linear regression model and regularizing them with sparsity-inducing penalties, such as the $L_1$ norm. Recently, \cite{li2023beyond} proposed a {truncated} SSVD optimization model that enforces sparsity and orthonormality using $L_0$ norm constraints, restricting the size of selected submatrices to be less than or equal to certain thresholds. The authors demonstrated that their SSVD approach generalizes some optimization models, including the well-known sparse PCA \citep{dey2023solving, del2023sparse}. To solve the SSVD problem, they developed exact and approximation algorithms. {However, as pointed out by the authors, SSVD requires careful selection of sparsity thresholds to yield meaningful results. More specifically, if these thresholds are not chosen appropriately, clustering may occur separately on rows and columns, failing to account for row-column interactions that characterize biclustering.}


Among the techniques available in the literature, the bipartite spectral graph partition approach and its variations have gained notable recognition \citep{dhillon2001co, kluger2003spectral, song2021weighted}.
Biclustering is related to bipartite graph partitioning, which depicts the relations between samples and features as a bipartite graph and translates the biclustering task as bipartite graph partitioning problem \citep{dhillon2001co, zha2001bipartite, rege2008bipartite, fan2010integer, huang2020auto}. 
From a graph-theoretical perspective, the data matrix can be viewed as a weighted complete bipartite graph, where the vertex set of one partition is the set of rows (samples) and the vertex set of the other partition is the set of columns (features). A weight incident on a row-column pair reflects the affinity between the corresponding sample and feature. Therefore, the biclustering problem involves finding specific types of 
complete bipartite subgraphs called \textit{bicliques} \citep{ames2014guaranteed}. A biclique in a bipartite graph corresponds to a dense region in the associated matrix, thereby representing a bicluster.
Many biclique extraction algorithms aim to identify the maximum biclique, which is the largest biclique within a given graph.
In contrast to the widely studied maximum clique problem \citep{wu2015review, hosseinian2020lagrangian}, the maximum biclique problem on bipartite graphs has two distinct variants: the maximum vertex biclique problem and the maximum edge biclique problem \citep{pandey2020maximum}. The maximum vertex biclique problem involves identifying a biclique with the highest number of vertices, and it can be efficiently solved in polynomial time \citep{garey1979computers}. On the other hand, the maximum edge biclique problem entails finding a biclique with the maximum number of edges and is proven to be NP-complete \citep{peeters2003maximum}.

The combinatorial nature of biclustering makes it computationally challenging, imposing severe restrictions on the practical use of global optimization approaches. To address such complexity, many mathematical programming formulations and heuristic algorithms have been proposed in the literature. In \cite{fan2010integer}, biclustering is studied via a mathematical optimization perspective, showing its connection with bipartite graphs and discrete optimization models. \cite{fan2010linear} propose a binary quadratic programming problem and three equivalent binary linear programming reformulations for some graph partitioning problems, including the bipartite one. \cite{fan2012multi} propose binary integer programming models to address the ratio cut \citep{hagen1992new} and normalized cut \citep{shi2000normalized} variants of the graph partitioning problem which are both NP-hard to solve. To overcome the computational complexity, they derive semidefinite programming (SDP) relaxations for both cuts. The authors extend these methods to bipartite graph partitioning models, specifically incorporating ratio cut and normalized cut measures. Moreover, they formulate the models as quadratically constrained programs to obtain globally optimal solutions. Computational results show that their exact method requires a significant amount of time to reach the optimal solution, as evidenced by experiments on graphs with 10 vertices.
As a model problem for biclustering, \cite{ames2014guaranteed} focuses on partitioning a
bipartite graph into dense disjoint subgraphs. This gives rise to the so-called \textit{$k$-densest-disjoint biclique} ($k$-DDB) problem which is NP-hard. Here, the objective is to identify $k$ disjoint bipartite complete subgraphs of the input graph such that the sum of their densities is maximized.
To address the $k$-DDB problem, a rank-constrained formulation is proposed and then a relaxation technique is applied, transforming it into an SDP program through matrix lifting. This relaxation approach is similar to that employed for traditional clustering problems, such as the graph partitioning problem \citep{wolkowicz1999semidefinite, sotirov2014efficient} and those proposed in recent papers focusing on the minimum sum-of-squares clustering problem \citep{piccialli2022sos, piccialli2022exact, piccialli2023global}. Ames shows that for certain instances, the globally optimal biclustering can be obtained from the optimal solution of the SDP relaxation. Specifically, the relaxation is tight when the edge weights in the input graph are concentrated on a set of disjoint bipartite subgraphs. {Due to the NP-hardness of the $k$-DDB problem, numerous heuristics have been proposed (see, e.g., the survey by \cite{metaheuristics} and the references therein). These methods use greedy and local search strategies, focusing on biclique enumeration. Among these, the SAMBA (Statistical-Algorithmic Method for Bicluster Analysis) algorithm, proposed by \cite{tanay2002discovering, tanay2004revealing}, is notably effective. SAMBA operates through several phases to achieve its results. In the preprocessing phase, the algorithm limits the number of vertices the bicliques can contain, thereby reducing the search space. It then performs an exhaustive enumeration on the restricted graph to find $k$ bicliques with the highest densities. Finally, in the postprocessing phase, the algorithm makes local improvements on the extracted bicliques by greedily adding or removing vertices. Applications of SAMBA on gene expression data show that it can discover biclusters in matrices containing approximately 5,000 genes and 500 conditions. However, it operates under bicluster size constraints and does not provide guarantees on the quality of the obtained solution.}



Although exact approaches are only applicable to small-scale instances, global optimization holds significance in bicluster analysis for two reasons. First, it helps to improve the accuracy of the results by ensuring that the biclusters found are truly representative of the patterns in the data. This is especially important when dealing with real-world applications, where local optimization techniques may not be sufficient to find the best biclusters. Second, a well-defined global optimization solver can be used to fine-tune existing heuristic approaches, identifying potential pitfalls, and thereby helping to develop novel approximate methods.

As a model problem for biclustering, we address the $k$-DDB problem and propose a global optimization algorithm based on the branch-and-cut technique. This tailored approach leads to a considerable increase of the problem size that can be solved, showing the solver's suitability for practical applications. The main contributions of the paper can be summarized as follows:


\begin{itemize}
    \item We design an SDP-based branch-and-cut algorithm for solving $k$-DDB instances to global optimality. For the upper bound routine we employ the SDP relaxation proposed by \cite{ames2014guaranteed} and strengthen it by adding valid inequalities. We use a first-order method to efficiently solve the resulting SDP relaxation in a cutting-plane fashion. 
    For the lower bound computation, we find feasible bicliques by designing a maximum weighted matching heuristic that exploits the solution of the SDP relaxation. Finally, we propose a tailored branching strategy where the subproblems can be reformulated as SDPs over lower dimensional positive semidefinite cones.
    \item Extensive numerical experiments conducted on both synthetic and real-world instances demonstrate that the proposed algorithm can handle graphs with a number of nodes approximately 20 times larger than those addressed by general-purpose exact solvers. To promote reproducibility and facilitate further research in the field, {the source code of the solver is publicly available at \url{https://github.com/INFORMSJoC/2024.0683} \citep{sudoso2024}}.
\end{itemize}

The remainder of the paper is organized as follows. Section \ref{sec:formulation} introduces key concepts from graph theory and illustrates the relationship between biclustering and the $k$-DDB problem. Section \ref{sec:bac_algorithm} outlines the components of the proposed branch-and-cut method. In particular, Section \ref{sec:sdp_relax} reviews the SDP relaxation and highlights its properties. Section \ref{sec:valid_ineq} presents a class of valid inequalities and the associated cutting-plane algorithm. Section \ref{sec:valid_ub} describes the upper bound computation using a first-order SDP solver and proves the validity of the obtained bounds. Section \ref{sec:heuristic} introduces the rounding heuristic. Section \ref{sec:branching} illustrates a branching strategy for reducing the size of child nodes. Section \ref{sec:comp_results} describes the implementation details and showcases computational results on both artificial and real-world instances. Finally, Section \ref{sec:conclusions} concludes the paper by discussing future research directions.

\section*{Notation}
Throughout the paper, $\mathbb{S}^n$ denotes the space of $n\times n$ real symmetric matrices, $\mathbb{R}^n$ is the space of $n$-dimensional real vectors and $\mathbb{R}^{m \times n}$ is the space of $m \times n$ real matrices.
We denote by ${X} \succeq 0$ a matrix ${X}$ that is positive semidefinite and by $\mathbb{S}^n_+$ be the set of positive semidefinite matrices of size $n \times n$. We denote by ${0_n}$ and ${1_n}$ the $n$-dimensional vectors of all zeros and of all ones, respectively. 
Given ${x} \in \mathbb{R}^n$, $\textrm{Diag}({x})$ is the $n \times n$ diagonal matrix with ${x}$ on its
diagonal, whereas given ${X} \in \mathbb{R}^{n \times n}$, $\textrm{diag}({X})$ is the vector with the diagonal elements of ${X}$. We denote by $\inprod{\cdot}{\cdot}$ the 
trace inner product. That is, for any
${A}, {B} \in \mathbb{R}^{m\times n}$, we define $\inprod{{A}}{{B}} := \textrm{tr}({B}^\top {A})$. Finally, given ${A} \in \mathbb{R}^{n \times n}$ we denote the set of its eigenvalues by $\lambda({A})$.

\section{Definitions and problem formulation}
\label{sec:formulation}
As mentioned above, data for biclustering is usually stored in a rectangular matrix with $n$ rows and $m$ columns. Let $A \in \mathbb{R}^{n \times m}$ represent such a matrix, where the rows are indexed by the set $S = \{1, \dots, n\}$ and the columns by the set $F = \{1, \dots, m\}$. Here, the entry $A_{ij}$ represents the relationship between row $i$ and column $j$. The definition of biclustering can be formalized as follows \citep{madeira2004biclustering, busygin2008biclustering}.

\begin{definition}
Given a data matrix $A \in \mathbb{R}^{n \times m}$, a biclustering of $A$ into $k \in \{2, \dots, \min\{n, m\}\}$ biclusters is a collection of row and column subsets $B = \{(S_1, F_1), \dots, (S_k, F_k)\}$ that satisfies specific characteristics of homogeneity and

(i) the collection $\{S_1, \dots, S_k\}$ forms a partition of the set of rows, i.e., 
\begin{align*}
    S_j \subseteq S \ : \ S_j \neq \emptyset, \ S_1 \cup \dots \cup S_k = S, \ S_j \cap S_{j'} = \emptyset \quad \forall j,j' \in \{1, \dots, k\}, \ j \neq j';
\end{align*}

(ii) the collection $\{F_1, \dots, F_k\}$ forms a partition of the set of columns, i.e., 
\begin{align*}
    F_j \subseteq F \ : \ F_j \neq \emptyset, \ F_1 \cup \dots \cup F_k = F, \ F_j \cap F_{j'} = \emptyset \quad \forall j,j' \in \{1, \dots, k\}, \ j \neq j'.
\end{align*}
\end{definition}

Therefore, a bicluster corresponds to a submatrix of the data matrix $A$. If $k$ is known, biclustering can be formulated as an optimization problem of the form $\max \{f(B): B \in \mathcal{P}(A, k)\}$, where $\mathcal{P}(A, k)$ is the set of all the collections of size $k$ containing row and column subsets of $A$, and the objective function $f : \mathcal{P}(A, k) \rightarrow \mathbb{R}$ is the merit function that define the type of homogeneity between rows and columns that we seek in each bicluster. 
Following the related literature, we cast this task as a partitioning problem on a bipartite graph. Before discussing the mathematical programming formulation, we define some concepts used in graph theory.
\begin{definition}
    Given a bipartite graph $G = ((U \cup V), E)$, a pair of disjoint independent subsets $U' \subseteq U$, $V' \subseteq V$ is a biclique of $G$ if the subgraph of $G$ induced by $(U' \cup V')$, denoted by $B' = G[(U' \cup V')]$, is complete bipartite. That is, $B'$ is a bipartite graph such that $(u, v) \in E$ for all $u \in U'$, $v \in V'$.
\end{definition}

Let $K_{n, m} =((U \cup V), E)$ denote a weighted complete bipartite graph, where the vertex sets are defined as $U = \{u_1, \dots, u_n\}$ and $V=\{v_1, \dots, v_m\}$. Here, the vertex $u_i \in U$ corresponds to the $i$-th row of $A$ and the vertex $v_{j} \in V$ corresponds to the $j$-th column of $A$. Let $w : U \times V \rightarrow \mathbb{R}$ be a weight function. There is an edge $(u_i, v_j) \in E$ with weight $w(u_i, v_j) = A_{ij}$ between each pair of vertices $u_i \in U$ and $v_j \in V$. Any biclique $(U' \cup V')$ of $K_{n,m}$ represents a \textit{row cluster} $U' \subseteq U$ and a \textit{column cluster} $V' \subseteq V$.
Given $k \in \{2, \dots, \min\{|U|, |V|\}\}$, the $k$-DDB problem aims to identify a set of $k$ disjoint bicliques $\{(U_1 \cup V_1), \dots, (U_k \cup V_k)\}$ of $K_{n, m}$ such that the sum of the densities of the complete subgraphs induced by these bicliques is maximized \citep{ames2014guaranteed}. Given the correspondence between biclustering and complete bipartite graphs, throughout the rest of the paper, we will use the terms ``biclusters'' and ``biclique'' interchangeably. The definition of density of a graph is given below.

\begin{definition}
    Given a weighted complete bipartite graph $K_{n, m} = ((U \cup V), E)$, the density of a subgraph $H = ((U' \cup V'), E')$ of $K_{n,m}$, denoted by $d_H$, is defined as the total edge weight incident at each vertex divided by the square root of the number of edges, i.e.,
\begin{align*}
    d_H = \frac{1}{\sqrt{|U'||V'|}} \sum_{u \in U', v \in V'} w(u, v).
\end{align*}

\end{definition}
The $k$-DDB problem can be formulated as a discrete optimization problem \citep{ames2014guaranteed}.  To this end, we consider the partition matrix ${X_U} \in \{0, 1\}^{n \times k}$ where the $i$-th row represents the cluster assignment for $u_i \in U$ and the $j$-th column is the characteristic vector of $U_j$, i.e., $({X_U})_{ij} = 1$ if $u_i \in U_j$ and 0 otherwise. Similarly, we consider the partition matrix ${X_V} \in \{0, 1\}^{m \times k}$ where the $i$-th row represents the cluster assignment for $v_i \in V$ and the $j$-th column is the characteristic vector of $V_j$, i.e., $({X_v})_{ij} = 1$ if $v_i \in V_j$ and 0 otherwise. 
Given a set of $k$ disjoint bicliques $\{(U_1 \cup V_1), \dots, (U_k \cup V_k)\}$, where $B_j = K_{n, m}[(U_j \cup V_j)]$ represents the complete subgraph of $K_{n, m}$ induced by the biclique $(U_j \cup V_j)$ for $j \in \{1, \dots, k\}$, the sum of their densities can be rewritten as
\begin{align*}
    \sum_{j=1}^k d_{B_j} &= \sum_{j=1}^k \frac{\sum_{u \in U_j, v \in V_j} w(u, v)}{\sqrt{|U_j||V_j|}} = \sum_{j=1}^k \frac{(X_U^\top {A} {X_V})_{jj}}{\sqrt{|U_j| |V_j|}}.
\end{align*} 
The resulting problem can be formulated as
\begin{subequations}
\label{prob:origin}
\begin{align}
\max \quad & \tr\left(({X_U^\top X_U})^{-\frac{1}{2}} X_U^\top {A} X_V ({X_V^\top X_V})^{-\frac{1}{2}}\right)  \\
\textrm{s.\,t.} \quad & {X_U} {1_k} = {1_n}, \ {X_V} {1_k} = {1_m}, \label{constr:row_sum}\\
& {X_U} \in \{0, 1\}^{n \times k}, \ {X_V} \in \{0, 1\}^{m \times k}.
\end{align}
\end{subequations}
Problem \eqref{prob:origin} is a quadratic fractional program with binary variables and linear constraints. The resulting partition matrices $X_U$ and $X_V$ are nonoverlapping and exclusive since constraints \eqref{constr:row_sum} ensure that each object in $U$ and $V$ is assigned to exactly one row-cluster and column-cluster, respectively. We now rewrite problem \eqref{prob:origin} as a non-convex quadratically constrained quadratic program (QCQP). Problem \eqref{prob:origin} can be reformulated by replacing the partition matrices $X_U$ and $X_V$ with the normalized partition matrices $Y_U$ and $Y_V$ as follows
\begin{subequations}
\label{prob:or}
\begin{align}
\max \quad & \tr(Y_U^\top A {Y_V})  \\
\textrm{s.\,t.} \quad & Y_U^\top {Y_U} = {I_k}, \ Y_V^\top {Y_V} = {I_k}, \label{constr:ort}\\
& Y_U Y_U^\top 1_n = 1_n, \ Y_V Y_V^\top 1_m = 1_m, \label{constr:ort_sum}\\
& {Y_U} \geq {0}_{n \times k}, \ {Y_V} \geq {0}_{m \times k}. \label{constr:nonneg}
\end{align}
\end{subequations}
\begin{proposition}
    Problems \eqref{prob:origin} and \eqref{prob:or} are equivalent.
\end{proposition}
\ifJOC
\proof{Proof.}
\else
\begin{proof}
\fi
We show that from any feasible solution of Problem \eqref{prob:origin} we can construct a feasible solution of Problem \eqref{prob:or} with the same objective value and vice versa. Let $X_U$ and $X_V$ be feasible Problem \eqref{prob:origin}. Denote by ${P_U} = ({X_U^\top X_U})^{-\frac{1}{2}} = \Diag(1/\sqrt{|U_1|}, \dots, 1/\sqrt{|U_k|})$ and ${P_V} = ({X_V^\top X_V})^{-\frac{1}{2}} = \Diag(1/\sqrt{|V_1|}, \dots, 1/\sqrt{|V_k|})$ the $k \times k$ matrices having on the diagonal the reciprocal of the square root of the sizes of $U$ and $V$ partitions, respectively. Thus, ${Y_U} = {X_U} {P_U}$ and ${Y_V} = {X_V} {P_V}$ represent the normalized partition matrices. One can easily verify that constraints \eqref{constr:ort}, \eqref{constr:ort_sum}, \eqref{constr:nonneg} hold by construction, and the objective function values coincide.
To prove the other direction, let $Y_U$ and $Y_V$ be feasible for problem \eqref{prob:or}. 
Constraints \eqref{constr:ort} imply that $Y_U$ and $Y_V$ are matrices with orthonormal columns. 
Using that $Y_U \geq 0$ entry wise, we obtain that
$(Y_U)_{ij} \neq 0$ implies that $(Y_U)_{ih} = 0$ for all $h \neq j$, so $Y_U$ has exactly one nonnegative entry per row. Similarly, using $Y_V \geq 0$ entry wise, we obtain that $(Y_V)_{ij} \neq 0$ implies that $(Y_V)_{ih} = 0$ for all $h \neq j$. The constraint $Y_U Y_U^\top 1_n = 1_n$ implies that the vector $1_n$ belongs to the span of the columns of $Y_U$. Therefore if $(Y_U)_{ij} \neq 0$ and $(Y_U)_{hj} \neq 0$ then $(Y_U)_{ij} = (Y_U)_{hj}$. Similarly, the constraint $Y_V Y_V^\top 1_m = 1_,$ implies that the vector $1_m$ belongs to the span of the columns of $Y_V$. 
Hence, constraints \eqref{constr:ort}, \eqref{constr:ort_sum}, \eqref{constr:nonneg} ensure that the columns of $Y_U$ and $Y_V$ can only assume two values: $(Y_U)_{ij} \in \left\{0, |U_j|^{-\frac{1}{2}}\right\}$ and $(Y_V)_{ij} \in \left\{0, |V_j|^{-\frac{1}{2}}\right\}$. 
Therefore, from $Y_U$ is always possible to construct a solution $X_U$ where $(X_U)_{ij} = 1$ if $(Y_U)_{ij} > 0$ and 0 otherwise; and from $Y_V$ a solution $X_V$ where $(X_V)_{ij} = 1$ if $(Y_V)_{ij} > 0$ and 0 otherwise. These matrices are feasible for \eqref{prob:origin} and achieve the same objective value of Problem \eqref{prob:or}.
\ifJOC
\hfill\Halmos
\endproof
\else
\end{proof}
\fi


Note that, while Problem \eqref{prob:origin} is a discrete optimization problem with non-convex objective function and linear constraints, Problem \eqref{prob:or} is a continuous formulation with both non-convex objective function and constraints. However, the latter exhibits a hidden discrete nature since combinatorial constraints are directly imposed on the continuous decision variables.
To find the globally optimal solution of Problem \eqref{prob:or}, one can use the spatial branch-and-bound algorithm implemented in general-purpose solvers such as Gurobi and BARON \citep{sahinidis1996baron}.
However, from a practical standpoint, non-convex QCQP problems are one of the most challenging optimization problems: the current size of instances that can be solved to provable optimality remains very small in comparison to other NP-hard problems \citep{d2003relaxations}. 
Moreover, even for small-scale problems (e.g., for graphs with $n + m = 50$ vertices), off-the-shelf solvers take a large amount of time to reach the globally optimal solution. Therefore, the objective of this paper is twofold: (i) to relax this intractable program to a tractable conic program that yields efficiently computable upper bounds on the original problem; (ii) to develop a tailored branch-and-cut algorithm that can solve large-scale instances that cannot be solved with standard software.


\section{Branch-and-cut algorithm}
\label{sec:bac_algorithm}
In this section, we present all the ingredients of the proposed exact algorithm. We start with an in-depth examination of the rank-constrained SDP formulation proposed by \citep{ames2014guaranteed} and describe its properties. This formulation will be useful in motivating the SDP relaxation, the definition of valid inequalities, the rounding heuristic, and the branching strategy.

\subsection{SDP relaxation}
\label{sec:sdp_relax}

The requirement for bipartite graph partitioning through the densest $k$-disjoint-biclique criterion is stated in the feasible set
\begin{align*}
\mathcal{P} = \left\{ \begin{bmatrix}
{Z_{UU}} & {Z_{UV}}\\
Z_{UV}^\top & {Z_{VV}}
\end{bmatrix} \in \mathbb{S}^{n+m} \ : \ \begin{array}{l}
     {Z_{    UU}} {1_n} = {1_n}, \ {Z_{VV}} {1_m} = {1_m}, \\
     \tr({Z_{UU}}) = k, \ \tr({Z_{VV}}) = k, \\
    {Z} \geq {0}, \ {Z} \in \mathbb{S}_+^{n+m} \\
    \end{array} \right\}.
\end{align*}
Thus, the equivalent low-rank SDP formulation proposed by \cite{ames2014guaranteed} is given by
\begin{equation}
\label{prob:SDPrank}
\begin{aligned}
\max_{{Z} \in \mathcal{P}} \quad & \frac{1}{2} \tr({W} {Z})\\
\textrm{s.t.} \quad & \rank({Z}) = k,
\end{aligned}
\end{equation}
where ${W} = \begin{bmatrix}
{0} & {A} \\
{A}^\top & {0}
\end{bmatrix} \in \mathbb{S}^{n+m}$ is the weighted adjacency matrix of $K_{n, m}$. 
The following proposition establishes the equivalence between the QCQP formulation and the rank-constrained SDP. {Although this result is presented in \cite{ames2014guaranteed}, we provide the formal proof here to offer a deeper understanding of the rank-constrained SDP formulation and its underlying structure.}
\begin{proposition}\label{th:equiv_bipartite}
    Problems \eqref{prob:or} and \eqref{prob:SDPrank} are equivalent.
\end{proposition}
\ifJOC
\proof{Proof.}
\else
\begin{proof}
\fi
We have to show that any feasible solution of Problem \eqref{prob:or} gives rise to a feasible solution of \eqref{prob:SDPrank} with
the same objective value and vice versa. We first transform the QCQP in \eqref{prob:or} into the rank-constrained SDP reformulation in \eqref{prob:SDPrank}. To this end, we lift ${Y_U}$ and ${Y_V}$ into the space of $(n+m) \times (n+m)$ matrices by introducing a symmetric block matrix
\begin{align*}
    {Z} = \begin{bmatrix}
{Z_{UU}} & {Z_{UV}}\\
Z_{UV}^\top & {Z_{VV}}
\end{bmatrix} \in \mathbb{S}^{n+m},
\end{align*}
where ${Z_{UU}} = {Y_U} Y_U^\top    \in \mathbb{S}^{n}$, $Z_{VV} = {Y_V} Y_V^\top \in \mathbb{S}^{m}$ and ${Z_{UV}} = {Y_U} Y_V^\top \in \mathbb{R}^{n \times m}$.
The objective function of Problem \eqref{prob:or} can be written as
\begin{align*}
\tr(Y_U^\top {A} {Y_V}) & = \frac{1}{2} \tr \left(
\begin{bmatrix}
Y_U^\top & {Y_V}^\top
\end{bmatrix}
\begin{bmatrix}
{0} & {A}\\
{A^\top} & {0}
\end{bmatrix}
\begin{bmatrix}
{Y_U}\\
{Y_V}
\end{bmatrix}
\right) = \frac{1}{2} \tr \left(
\begin{bmatrix}
{0} & {A}\\
{A^\top} & {0}
\end{bmatrix}
\begin{bmatrix}
{Y_U}\\
{Y_V}
\end{bmatrix}
\begin{bmatrix}
Y_U^\top & Y_V^\top
\end{bmatrix}
\right) \\
& = \frac{1}{2} \tr \left(
\begin{bmatrix}
{0} & {A}\\
{A^\top} & {0}
\end{bmatrix}
\begin{bmatrix}
Y_U Y_U^\top & Y_U Y_V^\top\\
Y_V Y_U^\top & Y_V Y_V^\top
\end{bmatrix}
\right) = \frac{1}{2} \tr({W}{Z}).
\end{align*}
Thus, both solutions attain the same objective value in their respective problems. Next, observe that since ${Z_{UU}}$ linearizes ${Y_U} Y_U^\top$ and ${Z_{VV}}$ linearizes ${Y_V} Y_V^\top$, constraints ${Z_{UU}}{1}_n = {1}_n$, $\tr({Z_{UU}}) = k$, ${Z_{VV}}{1_m} = {1_m}$ and $\tr({Z_{VV}}) = k$ hold by construction.
Moreover, ${Z} \geq {0}$ holds as well since ${Y_U}$ and ${Y_V}$ are elementwise nonnegative. 
{Let $\mathbbm{1}_{U_j} \in \{0, 1\}^n$ be the indicator vector of $U_j$, where the $i$-th component of $\mathbbm{1}_{U_j}$ is equal to 1 if $u_i \in U_j$ and 0 otherwise. Similarly, let $\mathbbm{1}_{V_j} \in \{0, 1\}^m$ be the indicator vector of $V_j$, where the $i$-th component of $\mathbbm{1}_{V_j}$ is equal to 1 if $v_i \in V_j$ and 0 otherwise.
After row and column reordering, $Z$ can be expressed as the sum of $k$ rank-1 matrices:}
\begin{align*}
    Z = 
\begin{bmatrix}
\sum_{j=1}^k \frac{1}{{|U_j|}} \mathbbm{1}_{U_j}  \mathbbm{1}_{U_j}^\top & \sum_{j=1}^k \frac{1}{\sqrt{|U_j||V_j|}} \mathbbm{1}_{U_j}  \mathbbm{1}_{V_j}^\top\\
\sum_{j=1}^k  \frac{1}{\sqrt{|U_j||V_j|}} \mathbbm{1}_{V_j}  \mathbbm{1}_{U_j}^\top & \sum_{j=1}^k \frac{1}{{|V_j|}} \mathbbm{1}_{V_j}  \mathbbm{1}_{V_j}^\top
\end{bmatrix}    
=
\sum_{j=1}^k 
\begin{pmatrix}
\frac{1}{\sqrt{|U_j|}} \mathbbm{1}_{U_j} \\
\frac{1}{\sqrt{|V_j|}} \mathbbm{1}_{V_j} 
\end{pmatrix}
\begin{pmatrix}
\frac{1}{\sqrt{|U_j|}} \mathbbm{1}_{U_j}\\
\frac{1}{\sqrt{|V_j|}} \mathbbm{1}_{V_j}
\end{pmatrix}^\top = \sum_{j=1}^k \omega_j \omega_j^\top.
\end{align*}
{Since the vectors $\omega_j$ are orthogonal to each other (each vertex in $U$ and $V$ is assigned to exactly one row-cluster and one column-cluster), they are linearly independent}. Therefore, $\rank(Z) = k$ and ${Z} \in \mathbb{S}_+^{n+m}$ hold.
To prove the other direction, we exploit the structure of any feasible solution ${Z}$ of Problem \eqref{prob:SDPrank}. One can easily verify that submatrices ${Z_{UU}}$ and ${Z_{VV}}$ encode pairwise row and column assignments, respectively. That is, the positivity of ${Z_{UU}}$ decides whether $u_i$ and $u_{i'}$ belong to the same row cluster: $({Z_{UU}})_{ii'} = 1/|U_j|$ if $u_i \in U_j$, $u_{i'} \in U_j$ and 0 otherwise. Analogously, the positivity of ${Z_{VV}}$ decides whether $v_i$ and $v_{i'}$ belong to the same column cluster: $({Z_{VV}})_{ii'} = 1/|V_j|$ if $v_i \in V_j$, $v_{i'} \in V_j$ and 0 otherwise. From this characterization, is possible to infer the corresponding matrices ${Y_U}$ and ${Y_V}$ in the following manner: if $u_i \in U_j$ then $({Y_U})_{ij} = 1/\sqrt{|U_j|}$ and 0 otherwise; if $v_i \in V_j$ then $({Y_V})_{ij} = 1/\sqrt{|V_j|}$ and 0 otherwise. Clearly, $Y_U^\top Y_U = I_k$, $Y_V^\top Y_V = I_k$, ${Y_U} Y_U^\top 1_n = 1_n$, ${Y_V} Y_V^\top 1_m = 1_m$ hold and the objective function values coincide. This concludes the proof.
\ifJOC
\hfill\Halmos
\endproof
\else
\end{proof}
\fi

Proposition \ref{th:equiv_bipartite} states that there exists a bijection among the feasible regions of Problems \eqref{prob:or} and \eqref{prob:SDPrank} that preserves the objective function value. This means that, from now on, we will focus on solving the rank-constrained SDP reformulation to global optimality.
The most interesting difference between Problem \eqref{prob:or} and Problem \eqref{prob:SDPrank} is that the nonlinear and non-convex objective function of \eqref{prob:or} becomes linear in \eqref{prob:SDPrank}. However, the feasible set of problem \eqref{prob:SDPrank} is still non-convex due to $\rank({Z}) = k$. 
By dropping the non-convex rank constraint, the resulting SDP relaxation  has a matrix variable of size $(n + m) \times (n + m)$ and $n+m+2$ linear constraints:
\begin{equation}
\label{prob:SDPrelax}
\begin{aligned}
\max_{{Z} \in \mathcal{P}} \quad & \frac{1}{2} \tr({W} {Z}).
\end{aligned}
\end{equation}
This kind of SDP is a double nonnegative program (DNN) since the block-matrix ${Z}$ is both positive semidefinite and elementwise nonnegative. 
Clearly, the optimal objective function value of the SDP relaxation provides an upper bound on the optimal value of the original combinatorial optimization problem. Furthermore, if the optimal solution of the relaxation has rank equal to $k$, then it is also optimal for problem \eqref{prob:SDPrank}. If this occurs, the relaxation not only provides a solution for the original problem, but also a certificate of its optimality thanks to the duality theory \citep{vandenberghe1996semidefinite}. Although there exist many polynomial time algorithms for solving problem \eqref{prob:SDPrelax}, general-purpose interior-point SDP solvers are not able to address large problems efficiently. A viable option for computing upper bounds for large graphs is further relaxing the SDP into an LP. Having already linear constraints and linear objective function, this can be done by replacing the constraint ${Z} \succeq {0}$ by ${Z} = {Z^\top}$. Unfortunately, the bound provided by the LP relaxation turns out to be weak when used in a branch-and-bound algorithm. A better option, instead, is to approximately solve problem \eqref{prob:SDPrelax} with an augmented Lagrangian method and then post-process the output of the SDP solver to get valid upper bounds. This approach will be discussed in Section \ref{sec:valid_ub}.

\subsection{Valid inequalities}
\label{sec:valid_ineq}
Problem \eqref{prob:SDPrelax} has been obtained by relaxing the non-convex rank constraint on ${Z}$, so the upper bound that it can provide may not be as tight as one could wish. Strong bounds play a key role in the branch-and-bound algorithm, as they allow for a more efficient exploration of the search space and a faster convergence to the optimal solution. As reported in the study by \cite{ames2014guaranteed}, the SDP relaxation is tight when the data is sampled from the so-called ``planted bicluster model''. Here, the data matrix ${A}$ is randomly generated in two steps. First, for each planted biclique $(U_h \cup V_h)$, the entries ${A}_{ij}$ are sampled independently and identically from a probability distribution $\Omega_1$ with mean $\omega_1$, for all $u_i \in U_h$ and $v_j \in V_h$. Next, the remaining entries of ${A}$ are sampled independently and identically from a distribution $\Omega_2$ with a mean $\omega_2$ that is much smaller than $\omega_1$. Under this model, recovery is guaranteed if the edge weights in the input graph are heavily concentrated on the edges of a few disjoint subgraphs. 
In real-world applications, the data matrix may not be compliant with the planted biclustering model, and thus the SDP relaxation may not be tight. To overcome this limitation, additional constraints can be imposed to strengthen the SDP bound.
In the following, we consider valid inequalities for the $k$-DDB problem. These are described in the work by \cite{de2022ratio} for the ratio-cut polytope and include the well-known pair and triangle inequalities associated with the cut polytope \citep{deza1997geometry}. Next, we discuss the separation routine within a cutting-plane algorithm.

\begin{proposition}\label{prop:ineq}
Fix a parameter $p \in \{1, \dots, k\}$. Then, the hypermetric inequalities
\begin{align}
\sum_{\substack{u_j, u_h \in Q_U \\ j < h}} (Z_{UU})_{jh} &\geq \sum_{u_j \in Q_U} (Z_{UU})_{ij} - (Z_{UU})_{ii} \quad \forall u_i \in U, \ \forall Q_U \subseteq U \setminus \{u_i\} : 1 \leq |Q_U| \leq p, \label{eq:ineq_U}\\
\sum_{\substack{v_j, v_h \in Q_V \\ j < h}} (Z_{VV})_{jh} &\geq \sum_{v_j \in Q_V} (Z_{VV})_{ij} - (Z_{VV})_{ii} \quad \forall v_i \in V, \ \forall Q_V \subseteq V \setminus \{v_i\} : 1 \leq |Q_V| \leq p \label{eq:ineq_V}
\end{align}
are valid for Problem \eqref{prob:SDPrank}.
\end{proposition}

\ifJOC
\proof{Proof.}
\else
\begin{proof}
\fi
    We first show the validity of inequality \eqref{eq:ineq_U}, following the proof structure in \cite{de2022ratio}. Denote by $\gamma$ the number of vertices in $Q_U$ that are in the same row cluster as $u_i$ and let $\Gamma$ be the size of the corresponding cluster. Clearly, $\gamma \geq 0$ and $\Gamma \geq 1$. Using the structure of Problem \eqref{prob:SDPrank} as described in Section \ref{sec:sdp_relax}, it can be verified that $\sum_{u_j \in Q_U} (Z_{UU})_{ij} = \frac{\gamma}{\Gamma}$, $(Z_{UU})_{ii} = \frac{1}{\Gamma}$ and $\sum_{{u_j, u_h \in Q_U : j < h}} (Z_{UU})_{jh} = \binom{\gamma}{2} \frac{1}{\Gamma}$ if $\gamma \geq 2$. Therefore, substituting these quantities in \eqref{eq:ineq_U}, we have
    \begin{align*}
        \begin{cases}
         \frac{\gamma}{\Gamma} \leq \frac{1}{\Gamma} & \text{if} \ \gamma \in \{0, 1\}, \\
        \frac{\gamma}{\Gamma} \leq \frac{1}{\Gamma} + \frac{\gamma (\gamma - 1)}{2\Gamma} & \text{otherwise}
         \end{cases}
    \end{align*}
    both of which are true. To prove the validity of inequality \eqref{eq:ineq_V} for the column clusters, it suffices to observe that $Z_{VV}$ shares the same structure as $Z_{UU}$ and, thus, by employing the same arguments, its validity holds as well.
\ifJOC
\hfill\Halmos
\endproof
\else
\end{proof}
\fi

Choosing $p = k$ {strengthens the SDP relaxation most among all $p \in \{1, \dots, k\}$}. However, it is important to note that there are $O(n^{p+1} + m^{p+1})$ inequalities of types \eqref{eq:ineq_U}-\eqref{eq:ineq_V} and hence, for large $p$, the SDP relaxation becomes too expensive to solve. Moreover, the separation problem over these inequalities is quite similar to the separation problem over clique inequalities, and thus it cannot be done efficiently unless $p$ is small. {For this reason, we focus on the case $p = 2$ where inequalities \eqref{eq:ineq_U}-\eqref{eq:ineq_V} capture all the pair and triangle inequalities described below}. 
Consider three vertices $u_i, u_j, u_h \in U$. If $(u_i, u_j)$ and $(u_j, u_h)$ are in the same row cluster $U' \subseteq U$, then $(u_i, u_h)$ must be in $U'$. This property holds for the column clusters as well. That is, given three vertices $v_i, v_j, v_h \in V$, if $(v_i, v_j)$ and $(v_j, v_h)$ are in the same column cluster $V' \subseteq V$, then $(v_i, v_h)$ must be in $V'$. 
Indeed, for $p=2$, we have:
\begin{align}
    ({Z_{UU}})_{ij} & \leq ({Z_{UU}})_{ii} &&\forall \ \textrm{distinct} \ u_i, u_j \in U\label{eq:pair_u},\\
    ({Z_{VV}})_{ij} & \leq ({Z_{VV}})_{ii} &&\forall \ \textrm{distinct} \ v_i, v_j \in V\label{eq:pair_v}, \\
({Z_{UU}})_{ij} + ({Z_{UU}})_{ih} & \leq ({Z_{UU}})_{ii} + ({Z_{UU}})_{jh} &&\forall \ \textrm{distinct} \ u_i, u_j, u_h \in U \label{eq:tri_u},\\
({Z_{VV}})_{ij} + ({Z_{VV}})_{ih} & \leq ({Z_{VV}})_{ii} + ({Z_{VV}})_{jh} &&\forall \ \textrm{distinct} \ v_i, v_j, v_h \in V.\label{eq:tri_v}
\end{align}
Here, the fact that two vertices $u_i$ and $u_j$ belong to the same row cluster is denoted by a strictly positive value for the entry $(Z_{UU})_{ij}$. If both $(Z_{UU})_{ij}$ and $(Z_{UU})_{ih}$ are strictly positive, they must both be equal to $(Z_{UU})_{ii}$; this forces $(Z_{UU})_{jh}$ to also be equal to $(Z_{UU})_{ii}$. Similarly, when two vertices $v_i$ and $v_j$ belong to the same column cluster we have a strictly positive value for the entry $(Z_{VV})_{ij}$. If both $(Z_{VV})_{ij}$ and $(Z_{VV})_{ih}$ are strictly positive, they must both be equal to $(Z_{VV})_{ii}$, thus forcing $(Z_{VV})_{jh}$ to be equal to $(Z_{VV})_{ii}$ as well.
We collect the above inequalities in the set $\textrm{MET} = \{Z \in \mathbb{S}^{n + m} : \textrm{inequalities} \ \eqref{eq:pair_u}{\textrm{-}}\eqref{eq:tri_v} \ \textrm{hold}\}$. Thus, the SDP relaxation becomes
\begin{equation}
\label{prob:SDPrelax_MET}
\begin{aligned}
\max_{} \quad & \frac{1}{2} \tr({W} {Z})\\
\textrm{s.t.} \quad & Z \in \mathcal{P} \cap \textrm{MET}.
\end{aligned}
\end{equation}
It is known that solving SDPs with a large number of inequalities is impractical for large graphs. Here, we propose to tackle this limitation via a cutting-plane algorithm that adds {hypermetric} inequalities only if they are violated. 
{Although the separation problem for pair and triangle inequalities can be solved in $O(n^3 + m^3)$ time through enumeration, an exact approach can be inefficient for large graphs. To mitigate this, the separation routine is performed heuristically, generating only a fixed number of violated inequalities during each cutting-plane iteration.}
In our cutting-plane approach, {hypermetric} inequalities are iteratively added and purged after each upper bound computation. First, the optimal solution of the basic SDP relaxation in \eqref{prob:SDPrelax} is computed. Next, an initial set of valid inequalities is added starting from the most violated ones. After the optimal solution is obtained, we purge all inactive constraints and new violated inequalities are added. The problem with an updated set of inequalities is solved and the process is repeated until the upper bound decreases. 
As we will show in the computational results, the SDP bounding routine combined with valid inequalities will be fundamental for closing the gap.

\subsection{Valid upper bounds}
\label{sec:valid_ub}

Using the optimal solution of an SDP relaxation within a branch-and-cut framework requires the computation of valid bounds. These bounds can be obtained by solving Problem \eqref{prob:SDPrelax_MET} to high accuracy.
Although interior-point methods can solve SDPs to arbitrary accuracy in polynomial time, they fail (in terms of time and memory requirement) when the scale of SDPs is large. On the other hand, first-order solvers based on augmented Lagrangian methods can scale to significantly larger problem sizes \citep{sun2015convergent, yang2015sdpnal}. However, it is well known that these methods can find solutions with moderate precision in a reasonable time, whereas progressing to higher precision can be very time-consuming \citep{oliveira2018admm}. In the following, we provide an efficient way to get a valid bound from the output of an augmented Lagrangian method that solves Problem \eqref{prob:SDPrelax_MET} to moderate accuracy. We adapt the method introduced by \cite{jansson2008rigorous} that adds a small perturbation to the dual objective function value.


Before showing how to compute valid upper bounds, we rewrite Problem \eqref{prob:SDPrelax_MET} as a semidefinite program in standard form which naturally leads to its dual formulation. Let $\mathcal{A}$ be linear map from $\mathbb{S}^{n+m}$ to $\mathbb{R}^{n+m+2}$ and $b \in \mathbb{R}^{n+m+2}$ defined as
\begin{equation*}
\mathcal{A}({Z}) \ : \  {Z} \to \begin{bmatrix}
           \inprod{Z_{UU}}{I_n} \\
           \frac{1}{2}({Z_{UU}} + Z_{UU}^\top) \\
            \inprod{Z_{VV}}{I_m} \\
           \frac{1}{2}({Z_{VV}} + Z_{VV}^\top)
         \end{bmatrix}, \quad b = \begin{bmatrix}
           k \\
           1_n \\
            k   \\
           1_m
         \end{bmatrix}.
\end{equation*}
Furthermore, we define the linear map $\mathcal{B} : \mathbb{S}^{n+m} \rightarrow \mathbb{R}^{q}$ to collect valid inequalities of the set MET, as described in Section \ref{sec:valid_ineq}.
The primal SDP can be written as
\begin{equation}
\label{prob:primal}
\begin{aligned}
\max \quad & \frac{1}{2} \tr({W} {Z}) \\
\textrm{s.t.} \quad & \mathcal{A}({Z}) = b, \ \mathcal{B}({Z}) \leq {0_q}, \\
& {Z} \geq {0}, \ {Z} \in \mathbb{S}^{n+m}_+.
\end{aligned}
\end{equation}
Consider Lagrange multipliers ${y_U} \in \mathbb{R}^n$, ${y_V} \in \mathbb{R}^m$, $\alpha_U, \alpha_V \in \mathbb{R}$, ${t} \in \mathbb{R}^q$, ${Q}$, ${S} \in \mathbb{S}^{n+m}$. Thus, using the standard derivation in Lagrangian duality theory, the dual of Problem \eqref{prob:primal} is given by
\begin{equation}
\label{prob:dual}
\begin{aligned}
\min \quad & {b^\top \lambda} = {y_U^\top} 1_n + y_V^\top 1_m + k(\alpha_U + \alpha_V) \\
\textrm{s.t.} \quad & \frac{1}{2} {W} - 
\mathcal{A}^\top({\lambda}) - \mathcal{B}^\top({t}) + {Q} + {S} = {0_{n \times m}}, \\
& {t} \geq {0_q}, \ {Q} \geq {0_{n \times m}}, \ {S} \in \mathbb{S}_+^{n+m},
\end{aligned}
\end{equation}
where ${\lambda} = [\alpha_U, {y_U}, \alpha_V, {y_V}] \in \mathbb{R}^{n+m+2}$ is the dual variable with respect to the affine constraints,
\begin{align*}
    \mathcal{A}^\top({\lambda}) = \begin{bmatrix}
           \frac{1}{2}({1_n} y_U^\top + {y_U} 1_{n}^\top) + \alpha_U I_n & {0_{n \times m}}\\
           {0_{m \times n}} & \frac{1}{2}({1_m} y_V^\top + {y_V} 1_m^\top) + \alpha_V I_m
         \end{bmatrix}
\end{align*}
is the adjoint operator of $\mathcal{A}$ and $\mathcal{B}^\top$ is the adjoint of $\mathcal{B}$. Before describing the post-processing we perform after each bound computation, we state a result bounding the eigenvalues of any feasible solution of Problem \eqref{prob:SDPrelax}.
\begin{proposition}\label{theorem:eigboundZ}
Let ${Z}$ be a feasible solution of Problem \eqref{prob:SDPrelax}, then $\lambda_{\max}({Z}) \leq 2$.
\end{proposition}
\ifJOC
\proof{Proof.}
\else
\begin{proof}
\fi
Let $\tilde{{Z}}_{UU} = 
\begin{bmatrix}
{Z_{UU}} & {0} \\
{0} & {0}
\end{bmatrix}$ and ${\tilde{Z}_{VV}} = 
\begin{bmatrix}
{0} & {0} \\
{0} & {Z_{VV}}
\end{bmatrix}$. From the decomposition lemma of positive semidefinite block matrices by \cite{bourin2012decomposition} we have 
\begin{equation*}
{Z} = 
{M} {\tilde{Z}_{UU}} {M^\top} +
{N} {\tilde{Z}_{VV}} {N^\top},
\end{equation*}
where ${M}, N \in \mathbb{R}^{(n+m) \times (n+m)}$ are orthogonal matrices. Therefore
\begin{align*}
\lambda_{\max}({Z}) = \max_{\|{v}\|=1} {v^\top} {Z} {v} &= \max_{\|{v}\|=1} {v^\top} \big({M} {\tilde{Z}_{UU}} {M^\top} +
{N} {\tilde{Z}_{VV}} {N^\top} \big) {v} \\
&\leq \max_{\|{v}\|=1} {v^\top} \big({M} {\tilde{Z}_{UU}} {M^\top} \big) {v} + \max_{\|{v}\|=1} {v^\top} \big(
{N} {\tilde{Z}_{VV}} {N^\top} \big) {v} \\
&= \lambda_{\max}\big( {M} {\tilde{Z}_{UU}} {M^\top} \big) + \lambda_{\max}\big( {N} {\tilde{Z}_{VV}} {N^\top}\big) \\
&= \lambda_{\max}\big({\tilde{Z}_{UU}} \big) + \lambda_{\max}\big( {\tilde{Z}_{VV}} \big) \\
&= \lambda_{\max}({Z_{UU}}) + \lambda_{\max}({Z_{VV}}) = 2,
\end{align*}
where the last equation holds since ${Z_{UU}}$ and ${Z_{UU}}$ are stochastic matrices, i.e., square matrices whose entries are nonnegative and whose rows all sum to 1.
\ifJOC
\hfill\Halmos
\endproof
\else
\end{proof}
\fi

The next lemma by \cite{jansson2008rigorous} is needed for proving the validity of the error bounds.
\begin{lemma}\label{lem:jansson}
Let ${S}, {X} \in \mathbb{S}^n$ be matrices that satisfy
$\lambda_{\min}({X}) \geq 0$ and $\lambda_{\max}({X}) \leq \bar{x}$ for some $\bar{x} \in \mathbb{R}$. 
Then the following inequality holds:
\begin{equation*}
    \inprod{{S}}{{X}} \geq \bar{x}\sum_{i \colon  \lambda_i({S}) <0}\lambda_i({S}).
\end{equation*}
\end{lemma}
We are now ready to state and prove the main result. Once the SDP has been solved approximately, the following theorem gives an upper bound on the optimal value of the SDP relaxation. This ensures the theoretical validity of the bounds produced within the branch-and-cut algorithm.
\begin{theorem}
\label{theorem:safe_bound}
Let ${Z^\star}$ be an optimal solution of \eqref{prob:primal} with objective function value $p^\star$. Consider the dual variables  ${y_U} \in \mathbb{R}^n$, ${y_V} \in \mathbb{R}^m$, $\alpha_U, \alpha_V \in \mathbb{R}$, ${t} \in \mathbb{R}^q$, ${Q} \in \mathbb{S}^{n+m}$, $Q \geq 0$. Set ${\tilde{S}} = -\frac{1}{2} {W} + \mathcal{A}^\top({\lambda}) + \mathcal{B}^\top({t}) - {Q}$. The following inequality holds:
\begin{equation*}
    p^\star \leq {y_U^\top} 1_n + y_V^\top 1_m + k(\alpha_U + \alpha_V) - 2 \sum_{i \colon  \lambda_i({\tilde{S}}) <0}\lambda_i({\tilde{S}}).
\end{equation*}
\end{theorem}

\ifJOC
\proof{Proof.}
\else
\begin{proof}
\fi
Let ${Z^\star}$ be optimal for the primal SDP \eqref{prob:primal}. Set ${\lambda} = [\alpha_U, {y_U}, \alpha_V, {y_V}]$, then
\begin{align*}
   \frac{1}{2}\inprod{{W}}{{Z^\star}} - {y_U^\top} 1_n - y_V^\top 1_m - k(\alpha_U + \alpha_V) - {t^\top} {0_q} &= \frac{1}{2}\inprod{{W}}{{Z^\star}} - \inprod{\mathcal{A}^\top({\lambda})}{{Z^\star}} - \inprod{\mathcal{B}^\top({t})}{{Z^\star}}\\
   & = \inprod{\frac{1}{2} {W} - \mathcal{A}^\top({\lambda}) - \mathcal{B}^\top({t})}{{Z^\star}} \\
   & = - \inprod{{Q} + {\tilde{S}}}{{Z^\star}}\\
   & = - \inprod{{Q}}{{Z^\star}} - \inprod{{\tilde{S}}}{{Z^\star}} \\ 
   & \leq - 2\sum_{i \colon  \lambda_i({\tilde{S}}) <0}\lambda_i({\tilde{S}}),
\end{align*}
where the last inequality holds since ${Q}$ is nonnegative and thanks to Lemma \ref{lem:jansson} where Proposition \ref{theorem:eigboundZ} with $\bar{x} = 2$ is used.
\ifJOC
\hfill\Halmos
\endproof
\else
\end{proof}
\fi

From a practical point of view, this post-processing step produces a safe overestimate for the optimal objective function value of the primal SDP in \eqref{prob:primal}, which is also an upper bound on the optimal value of the original rank-constrained problem in \eqref{prob:SDPrank}. More precisely, if the dual matrix ${\tilde{S}}$ is positive semidefinite (within machine accuracy), then $({y_U}, {y_V}, \alpha_U, \alpha_V, {t}, {\tilde{S}}, {Q})$ is a feasible solution of Problem \eqref{prob:dual} and ${y_U^\top} 1_n + y_V^\top 1_m + k(\alpha_U + \alpha_V)$ is already a valid upper bound. Otherwise, we increase the dual objective function value by adding the positive term $ - 2\sum_{i \colon  \lambda_i({\tilde{S}}) <0}\lambda_i({\tilde{S}})$ to it. In this way, we obtain a valid bound as proved by Theorem \ref{theorem:safe_bound}.
Note that for the bound computation, it is necessary to have an upper bound on the maximum eigenvalue of any feasible solution of Problem \eqref{prob:primal}. In our case, as shown in Proposition {\ref{theorem:eigboundZ}}, this value is known and is equal to 2 thanks to the special structure of Problem \eqref{prob:primal}.

\subsection{Heuristic}
\label{sec:heuristic}
In the previous sections, we discussed how to obtain strong upper bounds by means of semidefinite programming tools. A key step in solving a combinatorial optimization problem via a convex relaxation involves
rounding a solution of the relaxation to a solution of the original discrete
optimization problem. In this section, we develop a rounding algorithm that recovers a feasible biclustering, and thus a lower bound on the optimal value of Problem \eqref{prob:or}, from the solution of the SDP relaxation solved at each node. We emphasize that the rounding step is unnecessary when the relaxation admits a solution that is feasible for Problem \eqref{prob:SDPrank}. This phenomenon is known in the literature as ``exact recovery'' or tightness of the SDP relaxation. 
Indeed, if the SDP relaxation admits a rank-$k$ solution $Z^\star$, one can use the positive entries of $Z^\star_{UU}$ and $Z^\star_{VV}$ to recover the optimal row and column clusters and the positive entries of $Z^\star_{UV}$ to recover the optimal matching between them. As previously mentioned, the tightness of the SDP relaxation can be proven for input graphs that are compliant with the planted biclustering model \citep{ames2014guaranteed}. Unfortunately, for real-world problems, the SDP relaxation is not tight in general and hence a rounding procedure is needed. However, the solution of the SDP relaxation strengthened with valid inequalities is often a good approximation of the optimal $k$-DDB solution: as each cutting-plane iteration progresses, the optimal bicliques become more distinctly defined. This motivates exploiting the quality of the SDP solution for the design of the rounding heuristic. This procedure is shown in Algorithm \ref{alg:heuristic} and consists of the following steps.


\begin{algorithm}
\caption{Find feasible bicliques}
\label{alg:heuristic}
\KwIn{Data matrix $A$, number of biclusters $k$, solution $\bar{Z}$ of the SDP relaxation.}
\begin{enumerate}[leftmargin=*]
\item Cluster the rows of $\tilde{Z}_{UU}$ with $k$-means and record the rows partition $\tilde{X}_U$.
\item Cluster the rows of $\tilde{Z}_{VV}$ with $k$-means and record the columns partition $\tilde{X}_V$. 
\item Let $\tilde{W} = (\tilde{X}_U^\top \tilde{X}_U)^{-\frac{1}{2}}\tilde{X}_U^\top A \tilde{X}_V (\tilde{X}_V^\top \tilde{X}_V)^{-\frac{1}{2}}$. Solve the linear assignment problem
\begin{align*}
\centering
    \Delta^\star \in \argmax_{\Delta \in \{0, 1\}^{k \times k}} \left\{ \inprod{\Delta}{\tilde{W}} \ : \ \Delta 1_k = 1_k, \ \Delta^\top 1_k = 1_k \right\}.
\end{align*}
\item Reorder the columns of $\hat{X}_U$ according to $\Delta_{ij}^\star$ and set $U^\star_j \leftarrow \{u_i \in U \ : \ (\hat{X}_U)_{ij} = 1\}$ for all $j \in \{1, \dots, k\}$.
\item Reorder the columns of $\hat{X}_V$ according to $\Delta_{ij}^\star$ and set $V^\star_j \leftarrow \{v_i \in V \ : \ (\hat{X}_V)_{ij} = 1\}$ for all  $j \in \{1, \dots, k\}$.
\end{enumerate}
\KwOut{Bicliques $\{(U^\star_1 \cup V^\star_1), \dots, (U^\star_k \cup V^\star_k) \}$.}
\end{algorithm}

After solving the SDP relaxation at each node we obtain a solution $\bar{Z}$, which may not be feasible for the original low-rank problem. 
In Steps 1-2 we run the $k$-means algorithm \citep{lloyd1982least} on the rows of submatrices $\tilde{Z}_{UU}$ and $\tilde{Z}_{VV}$ to obtain the assignment matrices $\tilde{X}_U$ and $\tilde{X}_V$, where their columns represent the characteristic vectors of vertex sets $U$ and $V$, respectively. Let $\hat{U}_j = \{u_i \in U : (\hat{X}_U)_{ij} = 1\}$ and $\hat{V}_j = \{v_i \in V : (\hat{X}_V)_{ij} = 1\}$ for all $j \in \{1, \dots, k\}$ denote the obtained row and column clusters, respectively. Step 3 translates into an assignment problem, which is the problem of finding a matching in a weighted bipartite graph such that the sum of weights is maximized. Specifically, the goal is to assign each row cluster $\hat{U}_j$ to exactly one column cluster $\hat{V}_j$ such that the sum of the densities of the complete subgraphs of $K_{n, m}$ induced by bicliques $(\hat{U}_j \cup \hat{V}_j)$ is maximized. To this end, consider a weighted complete bipartite graph $\tilde{G} = ((\tilde{U} \cup \tilde{V}), \tilde{E})$ with vertex sets $\tilde{U} = \{\tilde{u}_1, \dots, \tilde{u}_k \}$ and $\tilde{V} = \{\tilde{v}_1, \dots, \tilde{v}_k \}$. 
Here, the vertex $\tilde{u}_i \in \tilde{U}$ corresponds the $i$-th row cluster $\hat{U}_i$, the vertex $\tilde{v}_j \in \tilde{V}$ corresponds the $j$-th column cluster $\hat{V}_j$, and the weight $\tilde{W}_{ij}$ of an edge $(\tilde{u}_i, \tilde{v}_j) \in \tilde{E}$ is the density of the complete bipartite subgraph of $K_{m,n}$ induced by $(\hat{U}_i \cup \hat{V}_j)$. We introduce an indicator variable $\Delta_{ij} \in \{0, 1\}$ for every edge $(\tilde{u}_i, \tilde{v}_j) \in \tilde{E}$, such that $\Delta_{ij} = 1$ if $(\tilde{u}_i, \tilde{v}_j)$ belongs to the perfect matching, or equivalently, if the $i$-th row cluster $\hat{U}_i$ and $j$-th column cluster $\hat{V}_j$ are assigned to the same bicluster, and $\Delta_{ij} = 0$ otherwise. 
To guarantee that the variables indeed represent a perfect matching, we add constraints to ensure that each vertex of $\tilde{G}$ is adjacent to exactly one edge in the matching. In other words, we ensure that each row and column clusters are assigned to exactly one bicluster.
Then, we can write the problem as an integer linear program (ILP). Note that even though the assignment problem in Step 3 constitutes an ILP, it can be solved in polynomial time because its constraint matrix is totally unimodular. Therefore, one can relax integrality constraints to $\Delta_{ij} \geq 0$ and restate the problem as an LP. In Steps 4-5 we reorder the columns of $\hat{X}_U$ and $\hat{X}_V$ according to the optimal assignment $\Delta^\star_{ij}$ and we output the corresponding bicliques $\{(U^\star_1 \cup V^\star_1), \dots, (U^\star_k \cup V^\star_k) \}$. 

\subsection{Branching subproblems}
\label{sec:branching}
Using cutting planes to tighten the SDP relaxation may not be sufficient to certify the optimality of the best biclustering solution found so far. If there are no violated cuts or the bound does not exhibit substantial improvement upon the inclusion of hypermetric inequalities, we terminate the cutting plane generation and proceed to branching. In other words, we partition the existing problem into smaller dimensional problems by fixing certain variables.
For the branching strategy, we exploit the structure of Problem \eqref{prob:SDPrank}. We make the following observations:
\begin{enumerate}
    \item If $u_i \in U$ and $u_j \in U$ are in the same row cluster, then $(Z_{UU})_{ih} = (Z_{UU})_{jh}$ for all $u_h \in U$ and $(Z_{UV})_{ih} = (Z_{UV})_{jh}$ for all $v_h \in V$. Similarly, if $v_i \in V$ and $v_j \in V$ are in the same column cluster, then $(Z_{VV})_{ih} = (Z_{VV})_{jh}$ for all $v_h \in V$ and $(Z_{UV}^\top)_{ih} = (Z_{UV}^\top)_{jh}$ for all $u_h \in U$. 
    \item If $u_i \in U$ and $u_j \in U$ are not in the same row cluster, then $(Z_{UU})_{ij} = 0$. Similarly, if $v_i \in V$ and $v_j \in V$ are not in the same column cluster, then $(Z_{VV})_{ij} = 0$.
\end{enumerate}
Accordingly, it is natural to produce a binary enumeration tree. Every time a node is split into two child nodes, a pair of vertices $(u_i, u_j) \subseteq U \times U$ or $(v_i, v_j) \subseteq V \times V$ is chosen and a cannot-link (CL) decision and a must-link (ML) decision are respectively imposed on the left and the right child. A CL decision over the two vertices forces them to be in different bicliques, while, on the other hand, an ML decision forces them to be in the same bicliques. Clearly, with these two types of branching decisions, the set of solutions associated with the parent node is partitioned into two disjoint subsets. This branching strategy draws inspiration from the one employed in exact solvers for one-way clustering methods \citep{du1999interior, piccialli2022sos}. 

Denote by $\textrm{ML}_U \subseteq U \times U$ the set of must-link decisions between vertices in $U$ and by $\textrm{ML}_V \subseteq V \times V$ be the set of must-link decisions between vertices in $V$. Let $c_U = |\textrm{ML}_U |$ and $c_V = |\textrm{ML}_V|$ be their cardinalities, respectively. Furthermore, let $\textrm{CL}_U \subseteq U \times U$ be the set of CL decisions between vertices in $U$ and $\textrm{CL}_V \subseteq V \times V$ be the set of CL decisions between vertices in $V$. The SDP relaxation that needs to be solved at any node of the tree is
\begin{subequations}
\label{prob:branching_or}
\begin{align}
\max_{} \quad & \frac{1}{2} \tr(W Z) \\
\textrm{s.t.} \quad 
& (Z_{UU})_{ih} = (Z_{UU})_{jh} \qquad \forall u_h \in U, \ \forall (u_i, u_j) \in \textrm{ML}_U, \label{constr:row_u1}\\
& (Z_{UV})_{ih} = (Z_{UV})_{jh} \qquad \forall v_h \in V, \ \forall (u_i, u_j) \in \textrm{ML}_U, \label{constr:row_u2}\\
& (Z_{VV})_{ih} = (Z_{VV})_{jh} \qquad \forall v_h \in V, \ \forall (v_i, v_j) \in \textrm{ML}_V, \label{constr:row_v1}\\
& (Z_{UV}^\top)_{ih} = (Z_{UV}^\top)_{jh} \qquad \forall u_h \in U, \ \forall (v_i, v_j) \in \textrm{ML}_V, \label{constr:row_v2}\\
& (Z_{UU})_{ij} = 0 \qquad \forall (u_i, u_j) \in \textrm{CL}_U, \\
& (Z_{VV})_{ij} = 0 \qquad \forall (v_i, v_j) \in \textrm{CL}_V, \\
& {Z} = \begin{bmatrix}
{Z_{UU}} & {Z_{UV}}\\
Z_{UV}^\top & {Z_{VV}}
\end{bmatrix}, \ Z \in \mathcal{P}.
\end{align}
\end{subequations}

Note that, an ML branching decision not only implies that $(Z_{UU})_{ij} > 0$ or $(Z_{VV})_{ij} > 0$, but
it also implies that two entire rows (and columns) are exactly the same. This suggests that, instead of adding $n + m$ additional constraints to force the $i$-th and $j$-th row of submatrices $Z_{UU}$ and $Z_{UV}$ or $Z_{UV}^\top$ and $Z_{VV}$ to be identical, the overall matrix $Z$ can be replaced by a smaller one. Of course, this replacement has to be correctly handled by defining a slightly different feasible set and objective function. Hence, it is convenient to generalize this problem, defining a generic form for the SDP relaxation to be solved at each node of the tree. Denote by $\sigma_U(i) \subseteq U$ the set of vertices in $U$ grouped with vertex $u_i \in U$ and by $\sigma_V(i) \subseteq V$ the set of vertices in $V$ grouped with vertex $v_i \in V$ according to the branching decisions leading to the current node. Consider two binary matrices $T_U \in \{0, 1\}^{(n - c_U) \times n}$ and $T_V \in \{0, 1\}^{(m - c_V) \times m}$ defined as
\begin{align*}
    (T_{U})_{ij} = \begin{cases}
    1 & \text{if } u_j \in \sigma_U(i)\\
    0              & \text{otherwise}
\end{cases}, \qquad \quad (T_{V})_{ij} = \begin{cases}
    1 & \text{if } v_j \in \sigma_V(i)\\
    0              & \text{otherwise}
\end{cases}.
\end{align*}
The generic form of the SDP relaxation is given by
\begin{subequations}
\label{prob:branching_shr}
\begin{align} 
\max~ \quad & \inprod{T_U A T_V^\top}{\Zbar_{UV}} \\
\st~ \quad & \Zbar_{UU} T_U 1_n = 1_{n-c_U}, \ \Zbar_{VV} T_V 1_m = 1_{m-c_V}, \label{eq:SDPellb}\\
& \inprod{T_U T_U^\top}{\Zbar_{UU}} = k, \ \inprod{T_V T_V^\top}{\Zbar_{VV}} = k, \label{eq:SDPellk}\\
& (\Zbar_{UU})_{ij} = 0 \qquad \forall (i, j) \in \textrm{CL}_U, \\
& (\Zbar_{VV})_{ij} = 0 \qquad \forall (i, j) \in \textrm{CL}_V, \\ 
& 
    \Zbar = \begin{bmatrix}
\Zbar_{UU} & \Zbar_{UV}\\
\Zbar_{UV}^\top & \Zbar_{VV}
\end{bmatrix}
\in \mathbb{S}_+^{n-c_U+m-c_V}, \ \Zbar \geq 0.
\end{align}
\end{subequations}

\begin{theorem}
    Problems \eqref{prob:branching_or} and \eqref{prob:branching_shr} are equivalent.
\end{theorem}

\ifJOC
\proof{Proof.}
\else
\begin{proof}
\fi
We start by showing that any feasible solution $\Zbar$ for problem \eqref{prob:branching_shr} can be transformed into a feasible solution for \eqref{prob:branching_or} with the same objective function
value. We define this equivalent solution $Z \in \mathbb{S}^{n + m}$ as
\begin{align*}
    Z = \begin{bmatrix}
T_U & 0 \\
0 & T_V
\end{bmatrix}^\top \begin{bmatrix}
\Zbar_{UU} & \Zbar_{UV}\\
\Zbar_{UV}^\top & \Zbar_{VV}
\end{bmatrix} \begin{bmatrix}
T_U & 0 \\
0 & T_V
\end{bmatrix} = \begin{bmatrix}
T_U^\top \Zbar_{UU} T_U & T_U^\top \Zbar_{UV} T_V\\
T_V^\top \Zbar_{UV}^\top T_U & T_V^\top \Zbar_{VV} T_V
\end{bmatrix} 
\end{align*}
By computing these products, one can easily verify that $T_U$ has the effect of expanding $\Zbar_{UU}$ and $\Zbar_{UV}$ to $n \times n$ and $n \times m$ matrices by replicating their rows according to the branching decisions between vertices in $U$. Similarly, $\Zbar_{VV}$ and $\Zbar_{UV}^\top$ are expanded by $T_V$ to $m \times m$ and $m \times n$ matrices by replicating their rows according to the branching decisions between vertices in $V$. Thus, constraints \eqref{constr:row_u1}, \eqref{constr:row_u2}, \eqref{constr:row_v1} and \eqref{constr:row_v2} hold by construction.
Moreover, $Z \geq 0$ and $Z \succeq 0$ hold by construction as well. Next, we have that 
\begin{align*}
    \tr(Z_{UU}) &= \inprod{I_n}{T_U^\top \Zbar_{UU} T_u} = \inprod{T_U T_U^\top}{\Zbar_{UU}} = k, \\
    \tr(Z_{VV}) &= \inprod{I_m}{(T_V)^\top \Zbar_{VV} T_V} = \inprod{T_V T_V^\top}{\Zbar_{VV}} = k, \\
    Z_{UU}1_n &= T_U^\top \Zbar_{UU} T_U 1_n = T_U^\top 1_{n-c_U} = 1_n, \\
    Z_{VV}1_m &= T_V^\top \Zbar_{VV} T_V 1_m = T_V^\top 1_{n-c_V} = 1_m.
\end{align*}
Finally, the objective function value is
\begin{align*}
    \frac{1}{2} \inprod{\begin{bmatrix}
0 & A\\
A^\top & 0
\end{bmatrix} }{\begin{bmatrix}
Z_{UU} & Z_{UV}\\
Z_{UV}^\top & Z_{VV}
\end{bmatrix}} &= 
 \frac{1}{2} \inprod{\begin{bmatrix}
0 & A\\
A^\top & 0
\end{bmatrix} }{\begin{bmatrix}
T_U & 0 \\
0 & T_V
\end{bmatrix}^\top \begin{bmatrix}
\Zbar_{UU} & \Zbar_{UV}\\
\Zbar_{UV}^\top & \Zbar_{VV}
\end{bmatrix} \begin{bmatrix}
T_U & 0 \\
0 & T_V
\end{bmatrix}} \\
& = \frac{1}{2} \inprod{ \begin{bmatrix}
T_U & 0 \\
0 & T_V
\end{bmatrix}     \begin{bmatrix}
0 & A\\
A^\top & 0
\end{bmatrix}     \begin{bmatrix}
T_U & 0 \\
0 & T_V
\end{bmatrix}^\top}{\begin{bmatrix}
\Zbar_{UU} & \Zbar_{UV}\\
\Zbar_{UV}^\top & \Zbar_{VV}
\end{bmatrix}} \\
& = \inprod{T_U A T_V^\top}{\Zbar_{UV}}.
\end{align*}
It remains to show that for any feasible solution $Z$ of problem \eqref{prob:branching_or} {it} is possible to associate a feasible solution for \eqref{prob:branching_shr} with the same objective function value. Before proceeding with the proof, we define matrices $C_U = T_U T_U^\top$ and $C_V = T_V T_V^\top$. Now, assuming that $Z$ is a feasible solution for problem \eqref{prob:branching_or}, we define a new block matrix
$\Zbar \in \mathbb{S}^{n - c_U + m - c_V}$ as
\begin{align*}
    \Zbar = \begin{bmatrix}
C_U^{-1} T_U  & 0\\
0 & C_V^{-1} T_V 
\end{bmatrix}  \begin{bmatrix}
Z_{UU} & Z_{UV}\\
Z_{UV}^\top & Z_{VV}
\end{bmatrix}     \begin{bmatrix}
T_U^\top C_U^{-1}  & 0\\
0 & T_V^\top C_V^{-1} 
\end{bmatrix},
\end{align*}
which is nonnegative and positive semidefinite. Then, we have
\begin{align*}
    \inprod{T_U T_U^\top}{\Zbar_{UU}} & =  \inprod{C_U}{C_U^{-1} T_U Z_{UU} T_U^\top C_U^{-1}} = \inprod{T_U^\top C_U^{-1} C_U C_U^{-1} T_U}{ Z_{UU} } = \inprod{I_n}{Z_{UU}} = k, \\
    \inprod{T_V T_V^\top}{\Zbar_{VV}} & =  \inprod{C_V}{C_V^{-1} T_V Z_{VV} T_V^\top C_V^{-1}} = \inprod{T_V^\top C_V^{-1} C_V C_V^{-1} T_V}{ Z_{VV} } = \inprod{I_m}{Z_{VV}} = k,\\
    \Zbar_{UU} T_U 1_n & = C_U^{-1} T_U Z_{UU} T_U^\top C_U^{-1} T_U 1_n = C_U^{-1} T_U Z_{UU} T_U^\top 1_{n-c_U} = C_U^{-1} T_U 1_n = 1_{n-c_U}, \\
    \Zbar_{VV} T_V 1_m & = C_V^{-1} T_V Z_{VV} T_V^\top C_V^{-1} T_V 1_m = C_V^{-1} T_v Z_{VV} T_V^\top 1_{m-c_v} = C_V^{-1} T_V 1_m = 1_{m-c_V}.
\end{align*}
Finally, we have to verify that $Z$ and $\Zbar$ have the same objective function value. Let
\begin{align*}
   \bar{A} = \begin{bmatrix}
T_U  & 0\\
0 & T_V 
\end{bmatrix}     \begin{bmatrix}
0 & A\\
A^\top & 0
\end{bmatrix}     \begin{bmatrix}
T_U  & 0\\
0 & T_V 
\end{bmatrix}^\top.
\end{align*}
Then, we have
\begin{align*}
    \inprod{T_U A T_V^\top}{\Zbar_{UV}} & = \frac{1}{2} \inprod{\bar{A}}{\begin{bmatrix}
\Zbar_{UU} & \Zbar_{UV}\\
\Zbar_{UV}^\top & \Zbar_{VV}
\end{bmatrix}} \\
& = \frac{1}{2} \inprod{\bar{A}}{\begin{bmatrix}
C_U^{-1} T_U  & 0\\
0 & C_V^{-1} T_V 
\end{bmatrix} \begin{bmatrix}
Z_{UU} & Z_{UV}\\
Z_{UV}^\top & Z_{VV}
\end{bmatrix}  \begin{bmatrix}
T_U^\top C_U^{-1}  & 0\\
0 & T_V^\top C_V^{-1} 
\end{bmatrix}} \\
& = \frac{1}{2} \inprod{   \begin{bmatrix}
0 & A\\
A^\top & 0
\end{bmatrix}  }{ \begin{bmatrix}
Z_{UU} & Z_{UV}\\
Z_{UV}^\top & Z_{VV}
\end{bmatrix}} = \frac{1}{2} \inprod{W}{Z}.
\end{align*}
\ifJOC
\hfill\Halmos
\endproof
\else
\end{proof}
\fi
Even though CL branching decisions have been left out while discussing the equivalence
of models \eqref{prob:branching_or} and \eqref{prob:branching_shr}, we stress that including such equality constraints does not affect the proof of equivalence: it holds as long as the sets $\textrm{CL}_U$ and $\textrm{CL}_V$ are properly updated whenever the solution matrix shrinks due to the addition of a new ML branching decision. Finally, note that, at the root node $\textrm{ML}_U = \textrm{CL}_U = \emptyset$ and $\textrm{ML}_V = \textrm{CL}_V = \emptyset$ imply that $T_U = I_n$ and $T_V = I_m$, and hence problem \eqref{prob:branching_shr} is equivalent to the basic SDP relaxation in \eqref{prob:SDPrelax}.

It remains to choose the pair of indices to branch on. A property that a feasible solution $Z$ of problem \eqref{prob:SDPrank} should have is that for each pair of vertices $u_i \in U$ and $u_j \in U$, it should be either $(Z_{UU})_{ij} = 0$ or $(Z_{UU})_{ii} = (Z_{UU})_{ij}$, that is
$(Z_{UU})_{ij} \left((Z_{UU})_{ii} - (Z_{UU})_{ij}\right) = 0$. Similarly, for each pair of vertices $v_i \in V$ and $v_j \in V$, it should be either $(Z_{VV})_{ij} = 0$ or $(Z_{VV})_{ii} = (Z_{VV})_{ij}$, that is
$(Z_{VV})_{ij} \left((Z_{VV})_{ii} - (Z_{VV})_{ij}\right) = 0$.
Hence, it makes sense to consider the violations of these quantities in choosing a pair of vertices in $U$ or $V$ to branch on.
More precisely, we select either a pair of vertices $(u_i, u_j) \in U \times U$ or $(v_i, v_j) \in V \times V$ and we put them in the same bicliques or in a different one. We define the sets
\begin{align*}
    I_U(u_i, u_j) = \left\{\bar{b}_{ij} \in \mathbb{R} \ : \ \bar{b}_{ij} = b_{ij} \cdot |U|, \ b_{ij} = \min_{u_i, u_j \in U} \left\{(Z_{UU})_{ij}, (Z_{UU})_{ii}-(Z_{UU})_{ij}  \right\} \right\}, \\
    I_V(v_i, v_j) = \left\{\bar{b}_{ij} \in \mathbb{R} \ : \ \bar{b}_{ij} = b_{ij} \cdot |V|, \ b_{ij} = \min_{v_i, v_j \in V} \left\{(Z_{VV})_{ij}, (Z_{VV})_{ii}-(Z_{VV})_{ij}  \right\} \right\}.
\end{align*}
Note that, the quantities $b_{ij}$ in $I_U$ and $I_V$ have been multiplied by $|U|$ and $|V|$, respectively, to make them comparable. Thus, the branching pair is selected as
\begin{align*}
(i^\star, j^\star) \in \argmax_{\substack{%
        u_i, u_j \in U\\
        v_i, v_j \in V
      }} \left\{ I_U(u_i, u_j) \cup I_V(v_i, v_j) \right\}.
\end{align*}
In other words, we choose two indices $i^\star$ and $j^\star$ with the least tendency to assign the corresponding vertices to the same biclique, or to different ones.

The overall branch-and-cut algorithm is shown in Algorithm \ref{alg:bbpseudocode}. It converges to a global maximum with optimality tolerance $\varepsilon$, due to the binary branching that performs an implicit enumeration of all the possible bicliques. In practice, a small number of nodes is required thanks to the strength of the upper bounding routine and the effectiveness of the heuristic, as will be shown in the next section.

\begin{algorithm}
\caption{Branch-and-cut algorithm for biclustering}
\label{alg:bbpseudocode}

\textbf{Input}: Data matrix $A$, number of biclusters $k$, optimality tolerance $\varepsilon$.

\begin{enumerate}[label*=\arabic*.] 
    \item Let $P_0$ be the initial $k$-DDB problem and set $\mathcal{Q} = \{P_0\}$. 
    \item Set $X_U^\star = \textrm{null}$, $X_V^\star = \textrm{null}$ with objective function value $v^\star = -\infty$.
    \item While $\mathcal{Q}$ is not empty:
    \begin{enumerate}[label*=\arabic*., rightmargin=15pt]
        \item Select and remove problem $P$ from $\mathcal{Q}$.
        \item Solve the SDP relaxation to get an optimal solution $Z$. Apply the post-processing procedure to get a valid upper bound $UB$.
        \item If $v^\star>-\infty$ and  $(UB - v^\star)/UB < \varepsilon$, go to \textit{Step 3}.
        \item Search for hypermetric inequalities violated by $Z$. If any are found, add them to the current SDP relaxation and go to \textit{Step 3.2}.
        \item Run the heuristic algorithm to get the matrix assignments $X_U$, $X_V$ and a lower bound $LB$. If $LB > v^\star$ then set $v^\star \leftarrow LB$, $X_U^\star \leftarrow X_U$, $X_V^\star \leftarrow X_V$.
        \item Select the branching pair $(i,j)$ and partition problem $P$ into two subproblems. For each problem update $T_U$, $T_V$ and $\textrm{CL}_U$, $\textrm{CL}_V$ accordingly, add them to $\mathcal{Q}$ and go to \textit{Step 3}.
    \end{enumerate}
\end{enumerate}
\textbf{Output}: Optimal assignments $X_U^\star$, $X_V^\star$ with objective value $v^\star$.
\end{algorithm}

\section{Computational results}
\label{sec:comp_results}
In this section, we first describe the implementation details and experimental setup. Next, we show and discuss the computational results on both artificial and real-world instances.

\subsection{Implementation details}
{\tt BICL-SDP}, which stands for Biclustering/Biclique SDP, is implemented in C++ with some routines written in MATLAB (version 2021b). The SDP relaxation at each node is solved by means of SDPNAL+, a MATLAB software that implements an augmented Lagrangian method for SDPs with bound constraints \citep{sun2020sdpnal+}. We set the accuracy tolerance of SDPNAL+ to $10^{-4}$ in the relative KKT residual. We also use Gurobi (version 10.0.2) to solve the LPs within the rounding heuristic. We run the experiments on a laptop with Intel(R) i7-13700H CPU, 32 GB of RAM, and Ubuntu 22.04.3 LTS. As for the cutting-plane setting, we randomly separate at most $100,000$ valid cuts, sort them in decreasing order with respect to the violation, and add the first $10,000$ violated ones in each iteration.
{We stop the cutting-plane procedure when the separation routine does not find violated inequalities or when the relative difference of the upper bound between two consecutive iterations is less than or equal to $10^{-3}$}. Moreover, the obtained set of inequalities at every node are inherited by the child nodes when branching is performed. Finally, we visit the tree with the best-first search strategy and require an optimality tolerance of {$\varepsilon = 0.001$} on the gap, i.e., we terminate the branch-and-cut method when {$\frac{UB - LB}{UB} < \varepsilon$}, where UB and LB denote the best upper and lower bounds, respectively. 


\subsection{Experiments on artificial instances}
As previously discussed, general-purpose solvers implement spatial branch-and-bound algorithms designed to tackle non-convex QCQPs. Therefore, the initial set of experiments aims to evaluate the scalability of state-of-the-art solvers when dealing with Problem \eqref{prob:or}. 
To this end, we compare {\tt BICL-SDP} against Gurobi on small-scale artificial graphs that are compliant with the planted biclustering model. That is, we generate input graphs where it is known a priori that $k$ disjoint bicliques with large edge weights exist. To create these graphs, we consider a weighted complete bipartite graph $K_{n, m}$ and select a set of $k$ disjoint bicliques $\{(U_1 \cup U_1), \dots, (U_k \cup V_k)\}$. Next, we construct a random matrix $A$ such that if $u_i \in U_h$, $v_j \in V_h$ for some $h \in \{1, \dots, k\}$, then $A_{ij}$ is sampled independently from the uniform distribution over the interval $[0, 1]$. If $u_i$ and $v_j$ belong to different bicliques of $K_{n, m}$, then $A_{ij} = 0$. Finally, we corrupt the entries of the obtained matrix by adding Gaussian noise with mean 0.5 and standard deviation $\sigma$. We generate random graphs with varying numbers of vertices $(n, m) \in \{10, 15, 20, 25\}$ for $n > m$, numbers of bicliques $k \in \{2, 3, 4\}$, and noise levels $\sigma \in \{0.1, 0.3\}$. Overall, the test set consists of 60 instances, labeled using the notation $\{n\}\_\{m\}\_\{k\}\_\{\sigma\}$. Note that in this scenario, the correct number of bicliques is known in advance, so we solve the instances only for that specific value of $k$. 

The experimental results are presented in Tables \ref{tab:artificial_results1} and \ref{tab:artificial_results3} for $\sigma = 0.1$ and $\sigma = 0.3$, respectively. The reported metrics include the computational time in seconds (Time [s]), the number of nodes (Nodes), and the percentage gap (Gap [\%]). We set a time limit of 3600 seconds and report the gap achieved along with the number of explored nodes upon reaching this limit. Additionally, for {\tt BICL-SDP}, we provide the number of cutting-plane iterations (CP) at the root node.

Table \ref{tab:artificial_results1} shows that Gurobi efficiently solves 20 out of 30 instances within the time limit. On the remaining instances, particularly those with a smaller total number of nodes and bicliques, it achieves small gaps. However, when $n+m \geq 40$ and $k = 4$, Gurobi reaches the time limit with a large gap.  In contrast, {\tt BICL-SDP} significantly outperforms Gurobi by solving all instances within 3 seconds at the root node. Notably, these instances are solved without cutting-plane iterations (CP = 0). This superiority is expected since, under the planted biclustering model and in scenarios without many outliers or numerous clusters, \cite{ames2014guaranteed} demonstrates that the SDP relaxation is tight with high probability. Turning to Table \ref{tab:artificial_results3}, we observe that Gurobi successfully solves 14 out of 30 instances and reaches the time limit for the remaining ones. Instances with $k=2$ exhibit smaller gaps, while larger graphs with $n+m \geq 40$ nodes and $k=4$ show notably larger gaps. In contrast, {\tt BICL-SDP} solves all instances within a few seconds. In 18 out of 30 instances, the SDP relaxation is tight at the root node (CP = 0). However, as the number of bicliques increases, the strength of the SDP relaxation diminishes, requiring some cutting-plane iterations and branching decisions for solving the instances. Specifically, on 5 instances, the cutting-plane algorithm at the root node is sufficient for closing the gap, with no need for branching (Nodes = 1). Conversely, on 7 instances, branching becomes necessary, with only a small number of nodes explored.

In summary, these findings complement those reported by \cite{ames2014guaranteed} on artificial graphs. Here, when the noise increases, the basic SDP relaxation fails to find and certify the globally optimal solution. Additionally, state-of-the-art global optimization solvers, such as Gurobi, demonstrate poor scalability as the number of nodes increases. Consequently, they cannot be effectively applied to instances encountered in practical applications, whose sizes far exceed those considered in artificial graphs. In the next section, we will extensively test the proposed {\tt BICL-SDP} on large-scale graphs and investigate how its behavior changes in real-world instances where information on the generative model and the correct value of $k$ are unavailable.

\begin{table}[!htbp]
    \centering
    \footnotesize
\begin{tabular*}{\textwidth}{@{\extracolsep{\fill}}lrrrrrrr@{\extracolsep{\fill}}}
\toprule%
& \multicolumn{3}{c}{Gurobi} & \multicolumn{4}{c}{BICL-SDP}\\
\cmidrule{2-4}\cmidrule{5-8}%
Instance & Time [s] & Nodes & Gap [\%] & CP & Time [s] & Nodes & Gap [\%]\\
\midrule
10\_10\_2\_0.1	&	3	&	401	&	< 0.1	&	0	&	2	&	1	&	< 0.1	\\
10\_10\_3\_0.1	&	58	&	5360	&	< 0.1	&	0	&	3	&	1	&	< 0.1	\\
10\_10\_4\_0.1	&	152	&	12431	&	< 0.1	&	0	&	2	&	1	&	< 0.1	\\
15\_10\_2\_0.1	&	7	&	562	&	< 0.1	&	0	&	2	&	1	&	< 0.1	\\
15\_10\_3\_0.1	&	104	&	4820	&	< 0.1	&	0	&	2	&	1	&	< 0.1	\\
15\_10\_4\_0.1	&	721	&	45812	&	< 0.1	&	0	&	2	&	1	&	< 0.1	\\
15\_15\_2\_0.1	&	31	&	2276	&	< 0.1	&	0	&	2	&	1	&	< 0.1	\\
15\_15\_3\_0.1	&	186	&	3123	&	< 0.1	&	0	&	2	&	1	&	< 0.1	\\
15\_15\_4\_0.1	&	3600	&	104814	&	0.13	&	0	&	2	&	1	&	< 0.1	\\
20\_10\_2\_0.1	&	13	&	931	&	< 0.1	&	0	&	2	&	1	&	< 0.1	\\
20\_10\_3\_0.1	&	298	&	8996	&	< 0.1	&	0	&	2	&	1	&	< 0.1	\\
20\_10\_4\_0.1	&	3600	&	103288	&	0.1	&	0	&	2	&	1	&	< 0.1	\\
20\_15\_2\_0.1	&	124	&	2475	&	< 0.1	&	0	&	2	&	1	&	< 0.1	\\
20\_15\_3\_0.1	&	809	&	13152	&	< 0.1	&	0	&	2	&	1	&	< 0.1	\\
20\_15\_4\_0.1	&	3600	&	18122	&	0.2	&	0	&	2	&	1	&	< 0.1	\\
20\_20\_2\_0.1	&	224	&	6417	&	< 0.1	&	0	&	2	&	1	&	< 0.1	\\
20\_20\_3\_0.1	&	1395	&	9885	&	< 0.1	&	0	&	2	&	1	&	< 0.1	\\
20\_20\_4\_0.1	&	3600	&	7372	&	7.71	&	0	&	3	&	1	&	< 0.1	\\
25\_10\_2\_0.1	&	132	&	2450	&	< 0.1	&	0	&	2	&	1	&	< 0.1	\\
25\_10\_3\_0.1	&	738	&	9885	&	< 0.1	&	0	&	2	&	1	&	< 0.1	\\
25\_10\_4\_0.1	&	3600	&	24416	&	0.19	&	0	&	2	&	1	&	< 0.1	\\
25\_15\_2\_0.1	&	201	&	8831	&	< 0.1	&	0	&	2	&	1	&	< 0.1	\\
25\_15\_3\_0.1	&	2768	&	19519	&	< 0.1	&	0	&	2	&	1	&	< 0.1	\\
25\_15\_4\_0.1	&	3600	&	5874	&	192.03	&	0	&	2	&	1	&	< 0.1	\\
25\_20\_2\_0.1	&	485	&	4806	&	< 0.1	&	0	&	2	&	1	&	< 0.1	\\
25\_20\_3\_0.1	&	3600	&	7847	&	154.38	&	0	&	2	&	1	&	< 0.1	\\
25\_20\_4\_0.1	&	3600	&	6306	&	194.75	&	0	&	2	&	1	&	< 0.1	\\
25\_25\_2\_0.1	&	1457	&	7752	&	< 0.1	&	0	&	2	&	1	&	< 0.1	\\
25\_25\_3\_0.1	&	3600	&	5753	&	144.03	&	0	&	2	&	1	&	< 0.1	\\
25\_25\_4\_0.1	&	3600	&	6993	&	23.19	&	0	&	2	&	1	&	< 0.1	\\
\bottomrule
\end{tabular*}
    \caption{Comparison between Gurobi and BICL-SDP on small-scale artificial instances with noise $\sigma = 0.1$. For each instance, we report the computational time in seconds (Time [s]), the number of nodes (Nodes), and the percentage gap (Gap [\%]). Additionally, for BICL-SDP, we include the number of cutting-plane iterations (CP) at the root node.}\label{tab:artificial_results1}
\end{table}

\begin{table}[!htbp]
    \centering
    \footnotesize
\begin{tabular*}{\textwidth}{@{\extracolsep{\fill}}lrrrrrrr@{\extracolsep{\fill}}}
\toprule%
& \multicolumn{3}{c}{Gurobi} & \multicolumn{4}{c}{BICL-SDP}\\
\cmidrule{2-4}\cmidrule{5-8}%
Instance & Time [s] & Nodes & Gap [\%] & CP & Time [s] & Nodes & Gap [\%]\\
\midrule
10\_10\_2\_0.3	&	9	&	1634	&	< 0.1	&	0	&	2	&	1	&	< 0.1	\\
10\_10\_3\_0.3	&	73	&	5145	&	< 0.1	&	0	&	2	&	1	&	< 0.1	\\
10\_10\_4\_0.3	&	3600	&	227576	&	2.19	&	1	&	16	&	13	&	< 0.1	\\
15\_10\_2\_0.3	&	35	&	6593	&	< 0.1	&	0	&	2	&	1	&	< 0.1	\\
15\_10\_3\_0.3	&	490	&	25430	&	< 0.1	&	1	&	3	&	1	&	< 0.1	\\
15\_10\_4\_0.3	&	3600	&	39510	&	4.52	&	1	&	17	&	15	&	< 0.1	\\
15\_15\_2\_0.3	&	90	&	2428	&	< 0.1	&	0	&	2	&	1	&	< 0.1	\\
15\_15\_3\_0.3	&	1151	&	52220	&	< 0.1	&	0	&	2	&	1	&	< 0.1	\\
15\_15\_4\_0.3	&	3600	&	13631	&	2.94	&	2	&	4	&	1	&	< 0.1	\\
20\_10\_2\_0.3	&	90	&	2457	&	< 0.1	&	0	&	2	&	1	&	< 0.1	\\
20\_10\_3\_0.3	&	3109	&	98938	&	< 0.1	&	1	&	3	&	1	&	< 0.1	\\
20\_10\_4\_0.3	&	3600	&	12406	&	4.17	&	2	&	14	&	7	&	< 0.1	\\
20\_15\_2\_0.3	&	486	&	8029	&	< 0.1	&	0	&	2	&	1	&	< 0.1	\\
20\_15\_3\_0.3	&	3600	&	37880	&	0.13	&	0	&	2	&	1	&	< 0.1	\\
20\_15\_4\_0.3	&	3600	&	5088	&	27.04	&	3	&	18	&	13	&	< 0.1	\\
20\_20\_2\_0.3	&	416	&	2594	&	< 0.1	&	0	&	2	&	1	&	< 0.1	\\
20\_20\_3\_0.3	&	3600	&	5499	&	127.53	&	0	&	2	&	1	&	< 0.1	\\
20\_20\_4\_0.3	&	3600	&	5901	&	212.64	&	1	&	3	&	1	&	< 0.1	\\
25\_10\_2\_0.3	&	353	&	7423	&	< 0.1	&	0	&	2	&	1	&	< 0.1	\\
25\_10\_3\_0.3	&	3600	&	48561	&	0.2	&	0	&	2	&	1	&	< 0.1	\\
25\_10\_4\_0.3	&	3600	&	5372	&	18.19	&	3	& 11	&	5	&	< 0.1	\\
25\_15\_2\_0.3	&	1354	&	7741	&	< 0.1	&	0	&	2	&	1	&	< 0.1	\\
25\_15\_3\_0.3	&	3600	&	7008	&	152.1	&	0	&	2	&	1	&	< 0.1	\\
25\_15\_4\_0.3	&	3600	&	5117	&	101.69	&	1	&	3	&	1	&	< 0.1	\\
25\_20\_2\_0.3	&	2251	&	8187	&	< 0.1	&	0	&	2	&	1	&	< 0.1	\\
25\_20\_3\_0.3	&	3600	&	4572	&	1.74	&	0	&	2	&	1	&	< 0.1	\\
25\_20\_4\_0.3	&	3600	&	5892	&	47.06	&	3	&	12	&	5	&	< 0.1	\\
25\_25\_2\_0.3	&	3496	&	10281	&	< 0.1	&	0	&	2	&	1	&	< 0.1	\\
25\_25\_3\_0.3	&	3600	&	5453	&	10.28	&	0	&	2	&	1	&	< 0.1	\\
25\_25\_4\_0.3	&	3600	&	6481	&	278.37	&	2	&	19	&	9	&	< 0.1	\\
\bottomrule
\end{tabular*}
    \caption{Comparison between Gurobi and BICL-SDP on small-scale artificial instances with noise $\sigma = 0.3$. For each instance, we report the computational time in seconds (Time [s]), the number of nodes (Nodes), and the percentage gap (Gap [\%]). Additionally, for BICL-SDP, we include the number of cutting-plane iterations (CP) at the root node.}
    \label{tab:artificial_results3}
\end{table}

\subsection{Experiments on real-world instances}
To better evaluate the effectiveness of {\tt BICL-SDP}, we consider several real-world instances from two sources: The Cancer Genome Atlas (TCGA) database \citep{weinstein2013cancer, cumbo2017tcga2bed, paci2017swim} and the datasets curated by \cite{de2008clustering}.
This publicly available benchmark collection comprises microarray datasets derived from experiments on cancer gene expression. In these gene expression datasets, most genes are lowly expressed and do not vary to a larger extent. Thus, all datasets have been pre-processed by removing non-informative genes (genes that do not display differential expression across samples). Finally, for each dataset, we build two expression matrices by sorting the genes by variances and selecting genes whose variances are in the second and third quartiles. The obtained expression matrices are denoted by suffixes ``v1'' and ``v2'', respectively. When applying biclustering algorithms to real-world data, sometimes there is no prior knowledge about the appropriate number of biclusters. In that case, a commonly employed method for selecting $k$ is the ``elbow method''. This approach is also widely used as a heuristic to determine the number of clusters in traditional one-way clustering methods. Here, we adopt this heuristic in the following manner. For each dataset, we solve the basic SDP relaxation at the root node for various values of $k \in [2, 10]$ and execute the rounding procedure to obtain a feasible biclustering solution. Next, we compute the objective function for different values of $k$ and generate a plot with the density of the corresponding $k$ bicliques on the y-axis and $k$ on the x-axis. We then identify the value of $k$ at the point of inflection, often referred to as the ``elbow'' of the curve. In cases where multiple inflection points are observed, we consider all the suggested values of $k$. Overall, we consider 64 instances with varying number of vertices ($(n+m) \in [248, 1248]$) and bicliques. Table \ref{tab:my_label} provides a summary of the cancer datasets, with $n$, $m$, and $k$ representing the number of genes, the number of samples, and the number of biclusters, respectively.

\begin{table}[!ht]
    \centering
    \footnotesize
    \begin{tabular}{l|lll}
    \toprule
        Instance &  $n$ & $m$ & $k$ \\
        \midrule
Bhattacharjee-2001-v1	&	386	&	203	&	4, 5	\\
Bhattacharjee-2001-v2 & 193 & 203 & 3, 4\\
Dyrskjot-2003-v1	&	601	&	40	&	2, 3	\\
Dyrskjot-2003-v2	&	301	&	40	&	2, 3, 4	\\
Golub-1999-v1	&	934	&	72	&	2	\\
Golub-1999-v2	&	467	&	72	&	2	\\
Gordon-2002-v1	&	813	&	181	&	2	\\
Gordon-2002-v2	&	407	&	181	&	2, 3	\\
Khan-2001-v1	&	534	&	83	&	2, 3	\\
Khan-2001-v2	&	267	&	83	&	2, 3	\\
Pomeroy-2002-v1	&	428	&	34	&	7, 8	\\
Pomeroy-2002-v2	&	214	&	34	&	7, 8	\\
Ramaswamy-2001-v1	&	681	&	190	&	2, 3, 4	\\
Ramaswamy-2001-v2	&	341	&	190	&	2, 3, 4\\
Risinger-2003-v1	&	885	&	42	&	2, 3	\\
Risinger-2003-v2	&	443	&	42	&	2, 3	\\
Singh-2002-v1	&	339	&	102	&	2, 4, 5	\\
Singh-2002-v2	&	169	&	102	&	2, 3	\\
Su-2001-v1	&	785	&	174	&	3, 4\\
Su-2001-v2	&	393	&	174	&	2, 3, 5	\\
\bottomrule
\end{tabular}\hspace{20pt}
    \begin{tabular}{l|lll}
    \toprule
        Instance &  $n$ & $m$ & $k$ \\
        \midrule
TCGA-BLCA-v1 & 822 & 426 & 2, 4\\
TCGA-BLCA-v2 & 411 & 426 & 2, 3\\
TCGA-BRCA-v1 & 942 & 206 & 2, 3\\
TCGA-BRCA-v2 & 471 & 206 & 2, 3\\
TCGA-CESC-v1 & 822 & 308 & 2 \\
TCGA-CESC-v2 & 411 & 308 & 2\\
TCGA-CHOL-v1 & 822 & 45 & 2, 3\\
TCGA-CHOL-v2 & 411 & 45 & 2, 3\\
TCGA-ESCA-v1 & 822 & 196 & 2, 3, 4\\
TCGA-ESCA-v2 & 411 & 196 & 2, 3, 4\\
TCGA-GBM-v1 & 822 & 172 & 2\\
TCGA-GBM-v2 & 411 & 172 & 2\\
\bottomrule
    \end{tabular}
    \caption{Cancer datasets used in the experiments. Instances from \citep{de2008clustering} are presented on the left, while those from TCGA \citep{weinstein2013cancer} are on the right. For each instance, we report the number of genes $n$, the number of samples $m$, and the number of target biclusters $k$.}
    \label{tab:my_label}
\end{table}

The experimental results are presented in Tables \ref{tab:results1} and \ref{tab:results2}. Here, we report the instance name (Instance),  the total number of vertices ($n+m$), the number of bicliques ($k$), the computational time in seconds (Time [s]), the number of branch-and-cut nodes (Nodes), and the final percentage gap (Gap [\%]). Moreover, we report some statistics relative to the root node. These include the percentage gap before performing cutting-plane iterations ($\textrm{Gap}_{\textrm{0}}$ [\%]), the number of cutting-plane iterations (CP), and the percentage gap after CP iterations ($\textrm{Gap}_{\textrm{CP}}$ [\%]).

\begin{table}[!ht]
\footnotesize
\centering
\begin{tabular}{lcc|rrr|rrr}
\toprule
Instance   & $n+m$ & $k$   &  $\textrm{Gap}_0$ [\%] & $\textrm{CP}$ & $\textrm{Gap}_{\textrm{CP}}$ [\%] & Nodes & Time [s] & Gap [\%]  \\
\midrule
Bhattacharjee-2001-v1 &	589	& 4	& 9.96	& 25 & 0.75 & 23 & 4627 & < 0.1\\
Bhattacharjee-2001-v1 &	589	& 5	& 20.41 & 26 & 0.99 & 31 & 6162 & < 0.1\\
Bhattacharjee-2001-v2 & 396 & 3 & 5.11 & 10 & 0.47 & 25 & 2099 & < 0.1\\
Bhattacharjee-2001-v2 & 396 & 4 & 10.99 & 11 & 0.58 & 17 & 1846 & < 0.1\\
Dyrskjot-2003-v1 & 641 & 2 & 7.79 & 15 & 0.64 & 15 & 8365 & < 0.1\\
Dyrskjot-2003-v1 & 641 & 3 & 6.34 & 28 & 0.58 & 13 & 9077 & < 0.1\\
Dyrskjot-2003-v2 & 341 & 2 & 6.24 & 8 & 0.87 & 9 & 909 & < 0.1\\
Dyrskjot-2003-v2 & 341 & 3 & 5.47 & 13 & 0.39 & 49 & 1753 & < 0.1\\
Dyrskjot-2003-v2 & 341 & 4 & 4.79	& 15 & 0.34 & 53 & 2492 & < 0.1\\
Golub-1999-v1 & 1006 & 2 & 8.90 & 46 & 0.45 & 5 & 20549 & < 0.1\\
Golub-1999-v2 & 539 & 2 & 7.03 & 18 & < 0.1 & 1 & 1696 & < 0.1\\   
Gordon-2002-v1 & 994 & 2 & 2.43 & 2 & 0.94 & 7 & 14304 & < 0.1\\
Gordon-2002-v2 & 588 & 2 & 2.03 & 3 & < 0.1 & 1 & 307 & < 0.1\\
Gordon-2002-v2 & 588 & 3 & 7.11 & 22 & 0.72 & 59 & 9559 & < 0.1\\
Khan-2001-v1 & 617 & 2  & 8.52 & 33 & 0.53 & 9 & 7884 & < 0.1\\
Khan-2001-v1 & 617 & 3  & 9.52 & 43 & < 0.1 & 1 & 4365 & < 0.1\\
Khan-2001-v2 & 350 & 2  & 6.85 & 8 & < 0.1 & 1 & 174 & < 0.1\\
Khan-2001-v2 & 350 & 3  & 9.52 & 13 & 0.19 & 7 & 611 & < 0.1\\
Pomeroy-2002-v1 & 462 & 7 & 10.74 & 13 & 0.54 & 11 & 2312 & < 0.1\\
Pomeroy-2002-v1 & 462 & 8 & 10.59 & 24 & 0.62 & 55 & 5579 & < 0.1\\
Pomeroy-2002-v2 & 248 & 7 & 10.06 & 6 & 0.47 & 9 & 431 & < 0.1\\
Pomeroy-2002-v2 & 248 & 8 & 9.62 & 7 & 0.57 & 39 & 1122 & < 0.1\\
Ramaswamy-2001-v1 & 871 & 2 & 7.63 & 27 & 0.27 & 5 & 9139 & < 0.1\\
Ramaswamy-2001-v1 & 871 & 3 & 7.19 & 31 & 0.23 & 3 & 7268 & < 0.1\\
Ramaswamy-2001-v1 & 871 & 4 & 7.76 & 27 & 0.71 & 35 & 14215 & < 0.1\\
Ramaswamy-2001-v2 & 531 & 2 &  7.15 & 14	& < 0.1 & 1 & 752 & < 0.1\\
Ramaswamy-2001-v2 & 531 & 3 &  6.90 & 12	& < 0.1 & 1 & 560 & < 0.1\\
Ramaswamy-2001-v2 & 531 & 4 &  7.92 & 19	& 0.25 & 21 & 3648 & < 0.1\\
Risinger-2003-v1 & 927 & 2 & 5.21 & 31 & < 0.1 & 1 & 5345 & < 0.1\\
Risinger-2003-v1 & 927 & 3 & 5.42 & 23 & 1.48 & 15 & 13879 & < 0.1\\
Risinger-2003-v2 & 485 & 2 & 5.01 & 7 & < 0.1 & 1 & 223 & < 0.1\\
Risinger-2003-v2 & 485 & 3 & 6.01 & 23 & 0.19 & 11 & 3014 & < 0.1\\
Singh-2002-v1 & 441 & 2 & 7.12 & 10 & < 0.1 & 1 & 482 & < 0.1\\
Singh-2002-v1 & 441 & 4 & 5.94 & 19 & 0.93 & 31 & 5351 & < 0.1\\
Singh-2002-v1 & 441 & 5 & 7.34 & 16 & 1.63 & 63 & 9968 & < 0.1\\
Singh-2002-v2 & 271 & 2 & 4.01 & 4 & < 0.1 & 1 & 73 & < 0.1\\
Singh-2002-v2 & 271 & 3 & 3.19 & 4 & < 0.1 & 1 & 75 & < 0.1\\
Su-2001-v1 & 959 & 3 & 5.59 & 10 & < 0.1 & 1 & 1623 & < 0.1\\
Su-2001-v1 & 959 & 4 & 3.98 & 17 & 0.25 & 5 & 10149 & < 0.1\\
Su-2001-v2 & 567 & 2 & 9.11 & 13 & 0.77 & 17 & 4456 & < 0.1\\
Su-2001-v2 & 567 & 3 & 5.61 & 11 & < 0.1 & 1 & 808 & < 0.1\\
Su-2001-v2 & 567 & 5 & 3.64 & 17 & < 0.1 & 1 & 1428 & < 0.1\\
\bottomrule
\end{tabular}
\caption{Computational results of {\tt BICL-SDP} on \cite{de2008clustering} instances. For each instance, we report the total number of vertices ($n+m$), the number of bicliques ($k$), the computational time in seconds (Time [s]), the number of branch-and-cut nodes (Nodes), and the final percentage gap (Gap [\%]). Moreover, we report following statistics at the root node: the percentage gap without cutting-plane iterations ($\textrm{Gap}_{\textrm{0}}$ [\%]), the number of cutting-plane iterations (CP), and the percentage gap after CP iterations ($\textrm{Gap}_{\textrm{CP}}$ [\%]).}
\label{tab:results1}
\end{table}

\begin{table}[!ht]
\footnotesize
\centering
\begin{tabular}{lcc|rrr|rrr}
\toprule
Instance   & $n+m$ & $k$   &  $\textrm{Gap}_0$ [\%] & $\textrm{CP}$ & $\textrm{Gap}_{\textrm{CP}}$ [\%] & Nodes & Time [s] & Gap [\%]  \\
\midrule
TCGA-BLCA-v1 & 1248 & 2 & 3.24	& 50 & 0.46 & 5 & 19626 & < 0.1\\
TCGA-BLCA-v1 & 1248 & 4 & 7.93	& 29 & 0.74 & 23 & 12774 & < 0.1\\
TCGA-BLCA-v2 & 837 & 2 & 2.14 & 19 & < 0.1 & 1 & 1440 & < 0.1\\
TCGA-BLCA-v2 & 837 & 3 & 3.35 & 20 & 0.54 & 19 & 4510 & < 0.1\\
TCGA-BRCA-v1 & 1148 & 2 & 0.82 & 16	& < 0.1 & 1 & 3798 & < 0.1\\
TCGA-BRCA-v1 & 1148 & 3 & 0.97 & 28 & 0.12 & 3 & 9521 & < 0.1\\
TCGA-BRCA-v2 & 677 & 2 & 1.67 &	7 & 0.31 & 3 & 1043 & < 0.1\\
TCGA-BRCA-v2 & 677 & 3 & 0.61 & 6 & < 0.1 & 1 & 358 & < 0.1\\
TCGA-CESC-v1 & 1130 & 2 & 2.01 & 31 & 0.52 & 13 & 18527 & < 0.1\\
TCGA-CESC-v2 & 719 & 2 & 2.09 & 18 & < 0.1 & 1 & 1271 & < 0.1\\
TCGA-CHOL-v1 & 867 & 2 & 0.89 & 8 & < 0.1 & 1 & 683 & < 0.1\\
TCGA-CHOL-v1 & 867 & 3 & 1.39 & 16 & 0.59 & 17 & 7080 & < 0.1\\
TCGA-CHOL-v2 & 456 & 2 & 0.58 & 1 & 0.1 & 1 & 22 & < 0.1\\
TCGA-CHOL-v2 & 456 & 3 & 1.32 & 14 & 0.23 & 7 & 921 & < 0.1\\
TCGA-ESCA-v1 & 1018 & 2 & 0.60	& 10 & < 0.1 & 1 & 949 & < 0.1\\
TCGA-ESCA-v1 & 1018 & 3 & 3.10 & 27 & 0.72 & 9 & 18405 & < 0.1\\
TCGA-ESCA-v1 & 1018 & 4 & 3.41 & 34 & 0.61 & 21 & 15619 & < 0.1\\
TCGA-ESCA-v2 & 607 & 2 & 0.39 & 3 & < 0.1 & 1 & 85 & < 0.1\\
TCGA-ESCA-v2 & 607 & 3 & 2.63 & 18 & 0.33 & 19 & 3566 & < 0.1\\
TCGA-ESCA-v2 & 607 & 4 & 2.20 & 15 & 0.28 & 25 & 2870 & < 0.1\\
TCGA-GBM-v1 & 994 & 2 & 3.63 & 33 & 0.54 & 3 & 13467 & < 0.1\\
TCGA-GBM-v2 & 583 & 2 & 2.30 & 22 & 0.14 & 5 & 2313 & < 0.1\\
\bottomrule
\end{tabular}
\caption{Computational results of {\tt BICL-SDP} on TCGA instances. For each instance, we report the total number of vertices ($n+m$), the number of bicliques ($k$), the computational time in seconds (Time [s]), the number of branch-and-cut nodes (Nodes), and the final percentage gap (Gap [\%]). Moreover, we report following statistics at the root node: the percentage gap without cutting-plane iterations ($\textrm{Gap}_{\textrm{0}}$ [\%]), the number of cutting-plane iterations (CP), and the percentage gap after CP iterations ($\textrm{Gap}_{\textrm{CP}}$ [\%]).}
\label{tab:results2}
\end{table}

Computational results show that {\tt BICL-SDP} successfully solves all instances within the required optimality gap (< 0.1\%). These findings confirm that adding valid inequalities to the SDP relaxation significantly reduces the gap. Specifically, 21 out of 64 instances are directly solved at the root node through cutting-plane iterations. For the remaining instances, a gap exists at the root node, and thus branching becomes necessary. Interestingly, only a few nodes are explored in these cases. Although the number of nodes in the branch-and-cut tree is small, the computational cost per node can be high due to the large number of inequalities. Looking at the statistics at the root node, larger instances require a higher number of cutting-plane iterations, resulting in generally longer computing times as the problem size increases. Notably, for most instances, the gap ($\textrm{Gap}_{\textrm{CP}}$ [\%]) becomes less than 1\% at the end of the cutting-plane procedure, often achieving a reduction of one order of magnitude compared to the initial gap ($\textrm{Gap}_{\textrm{0}}$ [\%]). As anticipated, the basic SDP relaxation is never tight at the root node, but the globally optimal solution is consistently found and certified by the proposed {\tt BICL-SDP}.
									
To better understand the contributions of both the cutting plane and the rounding heuristic at the root node, Figures \ref{fig:cp1} and \ref{fig:cp2} illustrate the percentage gap against the number of CP iterations. Note that, the rounding heuristic is executed after each upper bound computation to leverage the improved SDP solution. In the plots, the ``x'' marker denotes an update of the global lower bound at the corresponding CP iteration. These figures highlight two main findings. First, during the initial iterations, the gap diminishes at a faster rate. However, subsequent iterations remain valuable as they lead to further gap reduction. Second, the quality of the SDP bound significantly influences the rounding heuristic in generating feasible bicliques, as evidenced by the frequent updates of the lower bound. 
{Notably, in 57 out of 64 instances, the rounding heuristic finds the optimal solution early in the branch-and-bound process - at the root node and after just a few cutting plane iterations. This highlights the combined effectiveness of the cutting plane algorithm and the heuristic.}

Summarizing, {\tt BICL-SDP} can tackle instances having sizes approximately 20 times larger than those handled by off-the-shelf global optimization solvers. This achievement is attributed to several components of the proposed branch-and-cut algorithm: (i) the bounding procedure, which employs first-order methods to get computationally efficient SDP bounds; (ii) the fundamental role of the proposed inequalities within the cutting plane procedure; and (iii) the heuristic procedure complemented by the branching strategy which enables to find the global maximum either at the root node or within a limited number of explored nodes.

\begin{figure}[!ht]
    \centering
    \includegraphics[scale=0.64]{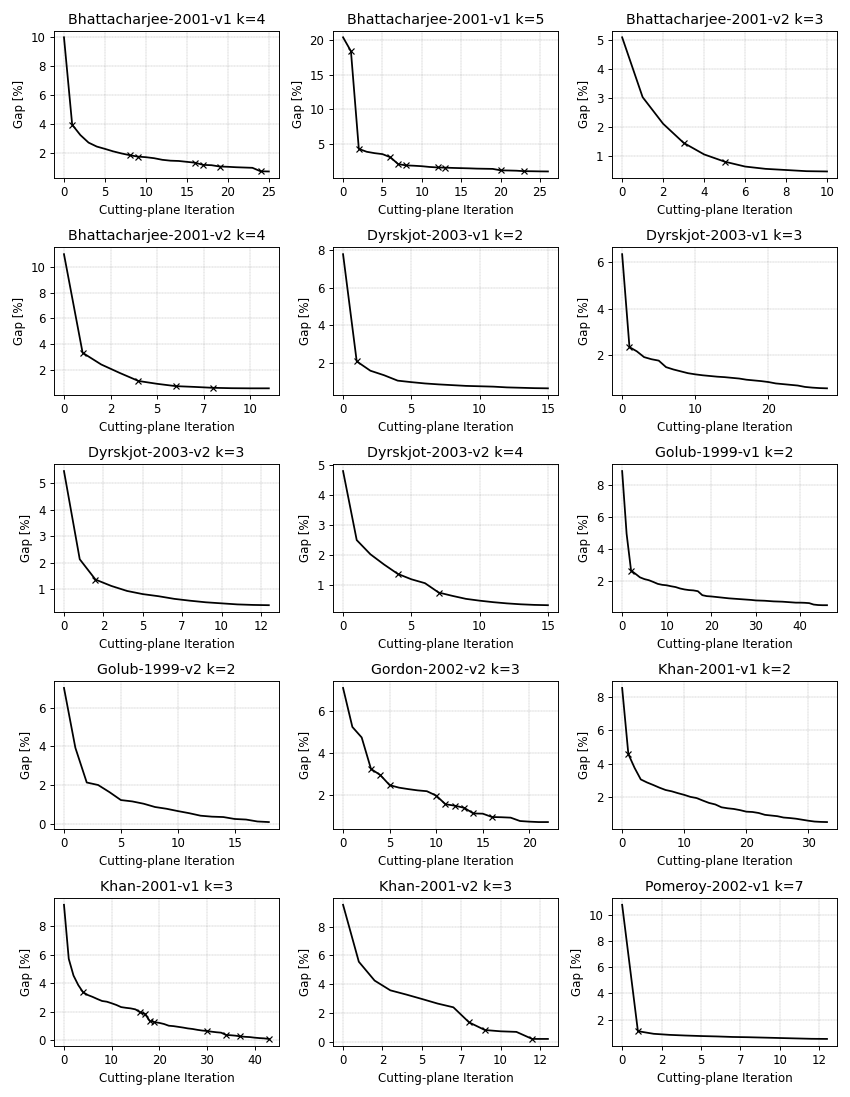}
    \caption{Percentage gap versus cutting-plane iterations at the root node. The “x” marker indicates that the global lower bound has been updated at the corresponding iteration.}
    \label{fig:cp1}
\end{figure}

\begin{figure}[!ht]
    \centering
    \includegraphics[scale=0.64]{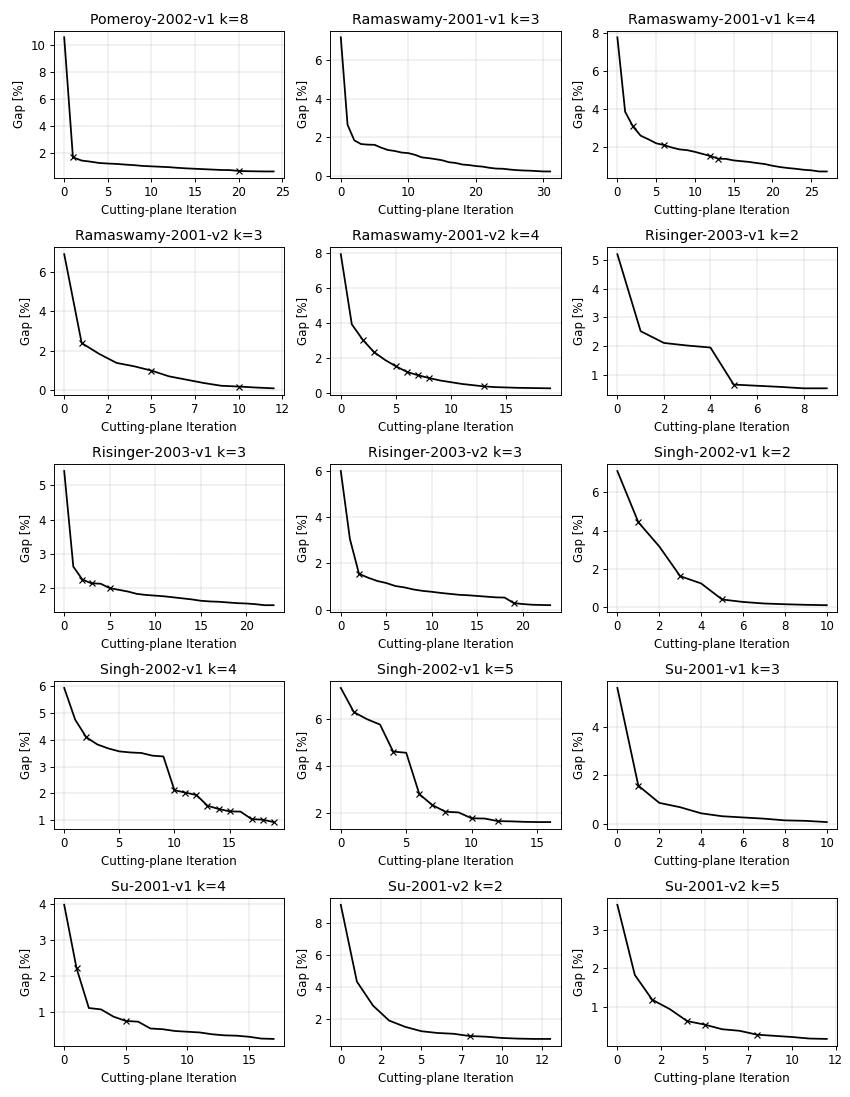}
    \caption{Percentage gap versus cutting-plane iterations at the root node. The “x” marker indicates that the global lower bound has been updated at the corresponding iteration.}
    \label{fig:cp2}
\end{figure}

\section{Conclusions}
\label{sec:conclusions}
This paper proposes a branch-and-cut algorithm for biclustering that exploits tools from semidefinite programming. 
Focusing on the $k$-DDB problem as a biclustering model, the objective is to identify a set of $k$ disjoint bicliques
within a given weighted complete bipartite graph such that the sum of the densities of the complete subgraphs induced by these bicliques is maximized. The upper bound routine employs first-order methods to get computationally efficient SDP bounds. To strengthen the bound, valid inequalities are added to the SDP relaxation in a cutting-plane fashion. Besides these upper bounds, the SDP relaxation also provides a primal solution that is used in the rounding heuristic to get feasible bicliques and hence lower bounds. Finally, a tailored branching strategy is used to reformulate subproblems as SDPs over lower dimensional positive semidefinite cones. Numerical results of the overall branch-and-cut algorithm impressively exhibit the efficiency of the proposed solver. To the best of the author's knowledge, no other exact solution method can handle real-world instances with {up to 1248 vertices}. This problem size is about 20 times larger than the one handled by general-purpose solvers. Clearly, as the number of vertices increases, the bottleneck shifts to solving the SDP relaxation with a large number of inequalities. Therefore, as a potential avenue for future research, one could develop a specialized SDP solver to further increase the size of biclustering instances that can be solved to provable optimality. {Additionally, future work could involve the exploration of LP-based post-processing techniques, as described in \cite{li2021strictly} and \cite{cerulli2021improving}, to further refine the quality of SDP bounds.}
Finally, the tools proposed in this paper could be extended to address variants of the biclustering problem, such as constrained biclustering \citep{pensa2008constrained, song2010constrained}. In constrained biclustering, user-defined constraints (e.g., the requirement for certain objects and/or attributes to be grouped together or kept separate) are incorporated into the problem. From an optimization perspective, this introduces several locally optimal solutions.

\section*{Acklowledgement}
I would like to thank Prof. Veronica Piccialli for her valuable suggestions and precious support.

\bibliography{biblio}

\end{document}